\documentclass[12pt]{amsart}
\usepackage{amsfonts, amsmath}
\usepackage{amssymb}
\usepackage{latexsym}
\usepackage{exscale}

\newtheorem{theo}{Theorem}[section]
\newtheorem{pro}[theo]{Proposition}
\newtheorem{lem}[theo]{Lemma}
\newtheorem{cor}[theo]{Corollary}

\newtheorem{rem}[theo]{Remark}

\newcommand{\calX}{{\mathcal{X}}}

\newcommand{\V}{{\mathcal{V}}}
\newcommand\av{\operatorname{av}}
\renewcommand\div{\operatorname{div}}
\newcommand\RR{\mathbb{R}}
\newcommand\CC{\mathbb{C}}
\newcommand\ZZ{\mathbb{Z}}

 \def\supp{\mathop{\rm supp}}
 \newcommand{\coinf}{C_{0}^\infty}
 \newcommand{\D}{{\mathcal{D}}}
 \newcommand{\schdot}{\dot H}
 \newcommand{\ext}{\widetilde H}
 \newcommand{\sch}{H}
 
 \newcommand{\ep}{\varepsilon}

\begin{document} 
\allowdisplaybreaks

\title[$L^p$ estimates for Schr\"odinger operators]{Maximal inequalities and Riesz transform estimates on $L^p$ spaces for Schr\"odinger operators with nonnegative potentials}

\author{Pascal Auscher}

\address{P. Auscher \& B. Ben Ali
\\
Universit\'e de Paris-Sud et CNRS UMR 8628
\\
91405 Orsay Cedex, France} \email{pascal.auscher@math.p-sud.fr}
 \email{besmath@yahoo.fr}

\author{Besma Ben Ali}


\date{\today}

\thanks{We thank B. Helffer for providing us with the unpublished reference \cite{HN} and also A. Ancona for indicating the relevance of \cite{IN}.}



\begin{abstract} We show various $L^p$ estimates for  Schr\"odinger operators  $-\Delta+V$  on $\RR^n$ and their square roots. We assume reverse H\"older estimates on the potential, and improve some results of Shen \cite{Sh1}.  Our main tools are  improved Fefferman-Phong inequalities and reverse H\"older estimates for weak solutions of $-\Delta+V$ and their gradients.

\end{abstract}

\maketitle

\begin{quote}
{\footnotesize\tableofcontents}
\end{quote}

\section{Introduction and main results}
Let $n\ge 1$ and $V$ be a locally integrable nonnegative function on $\RR^n$, not identically zero.
It is well-known that the following $L^1$ maximal inequality holds for $u \in \coinf(\RR^n)$,  real-valued,
\begin{equation}
\label{eq:L1}
  \|\Delta u \|_{1} + \| Vu \|_{1} \le  3\| - \Delta u + Vu \|_{1}.
\end{equation}
Here, $\|\ \|_{p}$ denotes the norm in $L^p(\RR^n)$. In fact, one has $\| Vu \|_{1} \le  \| - \Delta u + Vu \|_{1}$.

This follows either from work of Kato \cite{K2}, or  from work of Gallou\"et and Morel in semi-linear equations \cite{GM}. Nevertheless, we shall give a simple account of this.
This allows to define $-\Delta+V$ as an operator on $L^1(\RR^n)$ with domain $\D_{1}(\Delta) \cap \D_{1}(V)$. This was known before (\cite{V}) as we note that a similar inequality on $-\Delta+V+\lambda$ for some $\lambda\ge 0$  with constant depending also on $\lambda$ suffices. In this work we are interested in the possibility of having  ``homogeneous'' inequalities  ($\lambda=0$) or, equivalently,  on inequalities for $-\Delta+V+\lambda$ for $\lambda>0$ with constant independent of $\lambda$.

We turn to the $L^p$ theory for $1<p<\infty$. Assume that $V\in L^p_{loc}(\RR^n)$.  Then it is known that  $-\Delta+V$ a priori defined on $\coinf(\RR^n)$ is essentially $m$-accretive in  $L^p(\RR^n)$ (\cite{K1, K2, Se}) and the domain of the $m$-accretive extension contains  $\D_{p}(\Delta) \cap \D_{p}(V)=W^{2,p}(\RR^n)\cap L^p(\RR^n,V^p)$ as a dense subspace.    There are conditions to insure equality  in \cite{O, V, Da1, Si2}. But this is still not  enough to assert the validity of  the  $L^p$ version of  \eqref{eq:L1}, namely the \textit{a priori} inequality  for $u\in \coinf(\RR^n)$.
\begin{equation}
\label{eq:Lp}
  \|\Delta u\|_{p}+ \| Vu \|_{p} \lesssim \| - \Delta u + Vu \|_{p}.
\end{equation}
Here, $\sim$ is the equivalence in the sense of norms and $\lesssim$ the comparison of two norms.  
 A remark is that,  by standard Calder\'on-Zygmund theory, one can replace 
  $\|\Delta u \|_{p}$ by the equivalent quantity $ \|\nabla ^2 u\|_{p}$ as $1<p<\infty$. 

A natural question  is which condition on $V$ insures    \eqref{eq:Lp}.   An answer is the following. 
  
  \begin{theo}\label{th:1}  Let $1<q\le \infty$. If $V\in B_{q}$  then   for some $\ep>0$ depending only on $V$,  \eqref{eq:Lp} holds for $1<p<q+\ep$. \end{theo}

Here, $B_{q}$, $1<q\le\infty$,  is the class of the reverse H\"older weights: $w\in B_{q}$ if   $w \in L^q_{loc}(\RR^n)$,  $w>0$ almost everywhere and  there exists a constant $C$ such that  for all cube $Q$ of $\RR^n$,
\begin{equation}
\label{eq:bq}
\bigg( \frac 1 {|Q|} \int_{Q} w^q(x)\, dx\bigg)^{1/q} \le  \frac C {|Q|} \int_{Q} w(x)\, dx.
\end{equation}
If $q=\infty$, then the left hand side is the essential supremum on $Q$. 
The smallest $C$ is called the $B_{q}$ constant of $w$.  Examples of $B_{q}$ weights are the power weights $|x|^{-\alpha}$ for $-\infty<\alpha<n/q$ and positive polynomials for $q=\infty$. Note  that $B_{q}\subset B_{p}$ if $p<q$ and $w\in B_{q}$ implies $w\in B_{q+\ep}$ for some $\ep>0$ depending on the $B_{q}$ constant of $w$ (see \cite{Gra}). 

  Our result  extends the one of Shen obtained under the restriction that $n/2\le q$ and $n\ge 3$ \cite{Sh1}.  Prior to Shen's work, this was proved for positive polynomials when $p=2$  in \cite{HN}  and then when  $1<p<\infty$  in \cite{Gui1, Gui2,  Zh}.\footnote{After this work was completed, we learned of a new recent proof using representations via Lie groups  in \cite{DZ},  which  also  covers all positive fractional powers.}

A second family of inequalities concerns the square root (see below for definition). We recall at this point the identity
$$
\|\nabla  u\|_2^2 + \|V^{1/2} u\|_2^2 = \|(- \Delta  + V)^{1/2}u \|_2^2, \quad u\in \coinf(\RR^n).
$$
The a priori inequalities  
\begin{equation}
\label{eq:RT1}
\|\nabla u \|_{1,\infty} + \| V^{1/2}u \|_{1, \infty} \lesssim \| (- \Delta  + V)^{1/2}u \|_{1}
\end{equation}
 and 
 \begin{equation}
\label{eq:RTp}
\|\nabla u \|_{p} + \| V^{1/2}u \|_{p} \lesssim \| (- \Delta  + V)^{1/2}u \|_{p}
\end{equation} when  $1<p<2$ hold for $u\in C_{0}^\infty(\RR^n)$. Here, $\|\ \|_{p,\infty}$ is  the ``norm'' in the Lorentz space $L^{p,\infty}(\RR^n)$.
Actually, the first inequality  is attributed to Ouhabaz (unpublished) and the second one follows by interpolation. The proof of \eqref{eq:RT1} uses the fact that the heat kernel of $-\Delta+V$ is controlled pointwise by the one of $-\Delta$ and a theorem in \cite{DMc}. See \cite{DOY} where the needed estimates are proved and \cite{CD} where a similar argument is done for Riesz transforms on manifolds. See  also   \cite{Sik} for a different proof using finite speed of propagation for the wave equation.

We are interested in pushing the range of $p$ in \eqref{eq:RTp} beyond 2 and also in studying the converse inequalities, that is a priori  validity for smooth $u$  of
\begin{equation}
\label{eq:RRT1}
 \| (- \Delta  + V)^{1/2}u \|_{1,\infty} \lesssim  \|\nabla u \|_{1} + \| V^{1/2}u \|_{1}
\end{equation}
 and of
 \begin{equation}
\label{eq:RRTp}
 \| (- \Delta  + V)^{1/2}u \|_{p} \lesssim \|\nabla u \|_{p} + \| V^{1/2}u \|_{p}.
\end{equation}
Note that \eqref{eq:RTp} for $p$ implies \eqref{eq:RRTp} for the conjugate exponent $p'$. Hence, \eqref{eq:RRTp} already holds in the range $p>2$. 
The statement summarizing our results is the following.

\begin{theo}\label{th:2} 
\begin{enumerate}
\item  Let $V\in B_{q}$ for some $q>1$.  Then \eqref{eq:RTp} holds for $1<p<2(q+\ep)$.
   \item If $V\in A_{\infty}=\cup_{q>1}B_{q}$,  then \eqref{eq:RRT1} and \eqref{eq:RRTp} for $1<p<2$ hold.
   \item Let $V\in B_{q}$ for some $q>1$ and $q\ge n/2$.  Then  
   $\|\nabla u \|_{p}  \lesssim \| (- \Delta  + V)^{1/2}u \|_{p}$ holds
for  $1<p<q^*+\ep$   if $q<n$, and for $1<p<\infty$ if $q\ge n$.
\end{enumerate}
\end{theo}

Here,   $q^*= qn/(n-q)$ is the Sobolev exponent of $q$ if $q<n$. Note that $q^*\ge 2q$ exactly when $q\ge n/2$, hence  item 3 improves over item 1 for the gradient part. 
We note that    Shen  proved item 3  when  $n\ge 3$ and item 1 when $q\ge n/2$ and $n\ge 3$ \cite{Sh1}. 
 We shall fully prove this theorem, even item 3 with  an argument of a different nature that is interesting in its own right.
 
Note that one can also prove inequalities similar to \eqref{eq:RTp} for fractional powers 
$(-\Delta+V)^{-s}$, $0<s<1$, with range $1<p<(q+\ep)/s$. We shall not pursue this here.

Our results are satisfactory for reverse H\"older potentials as they make a bridge with
the known results for $L^1_{loc}$ nonnegative potentials. Let us list some other consequences to illustrate this. 

\begin{cor}\label{cor:1} Let $n\ge 1$, $1<p<\infty$ and $V\in B_{p}$.
Then  the $m$-accretive extension on $L^p(\RR^n)$ of  $-\Delta+V$ defined  on $C_{0}^\infty(\RR^n)$ has domain equal to  $\D_{p}(\Delta) \cap \D_{p}(V)$. In particular, for $p=2$, $-\Delta+V$ defined on $H^2(\RR^n)\cap L^2(\RR^n, V^2)$ is self-adjoint in $L^2(\RR^n)$.  
\end{cor}

  This applies to power weights $c|x|^{-\alpha}$  although this particular application is known by other methods \cite{O}.

\begin{cor}\label{cor:2} Let $n\ge 1$. Assume  $V\in A_{\infty}$ and $1<p<2$ or $V\in B_{p/2}$ and $2<p<\infty$, then $(-\Delta+V)^{1/2}$ has $L^p$-domain  equal to  $\D_{p}((-\Delta)^{1/2}) \cap \D_{p}(V^{1/2})= W^{1,p}(\RR^n)\cap L^{p}(\RR^n, V^{p/2})$. \end{cor}

Further  easy consequences are  the following estimates.
Set $\widehat p=\sup (2p,p^*)$ for $1<p<\infty$ with $p^*=\infty$ if $p\ge n$.

 \begin{cor} \label{cor:3} Assume that $V \in B_q$ for some $q>1$.  Then for $\ep>0$ depending only on $V$,
  \begin{enumerate}
  \item $V^{1/2} H^{-1} V^{1/2}$ is bounded on $L^p(\RR^n)$ for $(2(q+\ep))'<p<2(q+\ep)$.
  \item $V^{1/2} H^{-1} \div$  is bounded on $L^p(\RR^n)$ for $(\widehat{q+\ep})'<p<2(q+\ep)$.
  \item $\nabla  H^{-1} V^{1/2}$ is bounded on $L^p(\RR^n)$ for $(2(q+\ep))'<p< \widehat{q+\ep}$.
  \item $\nabla H^{-1} \div $ is bounded on $L^p(\RR^n)$ for $(\widehat{q+\ep})'<p< \widehat{q+\ep}$.
\end{enumerate}
 \end{cor}
 
 Again, this result extends  the ones of Shen in \cite{Sh1} obtained with the restriction $q\ge n/2$ and $n\ge 3$. He also proved  bounds for 
$V^{1/2}\nabla H^{-1}$, which we can recover  by our methods under the same hypotheses (and for $n\ge 1$ instead of $n\ge 3$). We therefore do not include such results. 
\

 We mention without proof that our results admit local versions, replacing $V\in B_{q}$ by  $V\in B_{q,loc}$ which is defined by the same conditions on cubes with sides less than 1. Then we get the corresponding results and estimates  for $H+1$ instead of $H$.  The results on operator domains are valid under local assumptions.

\

Our arguments are based on local estimates and this is fortunate because there is no auxiliary global weight as in \cite{Sh1}.   Our main  tools are

1) An improved Fefferman-Phong inequality for $A_{\infty}$ potentials.

2) Criteria for proving  $L^p$ boundedness of operators in absence of kernels. 

3)  Mean value inequalities for  nonnegative subharmonic functions against $A_{\infty}$ weights.

4) Complex  interpolation, together with  $L^p$ boundedness of imaginary powers 
of $-\Delta+V$ for $1<p<\infty$. 

5) A Calder\'on-Zygmund decomposition adapted to level sets of the maximal function of $|\nabla f| + |V^{1/2}f| $.

6) Reverse H\"older inequalities involving $\nabla u$ and $V^{1/2} u$ for weak solutions
of $-\Delta u + Vu=0$.  

\

The latter estimates are of independent interest and we give a rather complete picture. This is more than necessary
for applications to the inequality \eqref{eq:RTp}.

\section{An improved Fefferman-Phong inequality}\label{sec:FP}

Usual Fefferman-Phong inequalities take the form
$$
\int_{\RR^n} m(x)^2 |u(x)|^2\, dx \le  C \int_{\RR^n}    |\nabla u(x)|^2 + w(x)|u(x)|^2\, dx
$$
for $u\in C_{0}^\infty(\RR^n)$ where $m$ is a positive weight function depending on the potential $w$. If  $w \in B_{q}$ and $q\ge n/2$, there is such a function $m$  \cite{Sh1}. If $q<n/2$,  it is not clear how to define $m$ in function of $w$. Nevertheless, local inequalities on cubes $Q$ still hold  and depend on the scaling defined by the quantity $R^2 \av_{Q}w$ (The notation $\av_{E}v$ means $\frac{1}{|E|}\int_{E}v$).

\begin{lem}\label{lem:FP} Let $w\in A_\infty$ and $1\le p<\infty$. Then there are constants $C>0$ and $\beta\in (0,1)$ depending only
on the $A_\infty$ constant of $w$, $p$ and $n$ such that for all cubes $Q$ $($with sidelength $R)$ and $u \in C^1(\RR^n)$, one has
$$
\int_Q |\nabla u|^p + w|u|^p \ge \frac{C m_{\beta}({R^p \av_Q w})}{ R^{p}} \int_Q |u|^p 
$$
where $m_{\beta}(x)  = x$ for $x\le 1$ and $m_{\beta}(x) = x^\beta$ for $x\ge 1$.  
\end{lem}

This lemma with $\beta=0$ is already in \cite{Sh1} when $p=2$. The improvement occurs when $R^p \av_{Q}w\ge 1$ and  is crucial for us in Section \ref{sec:weaksol}. Such an improvement has also applications to  criteria for compactness of resolvents for magnetic Schr\"odinger operators (personal communication of B. Helffer).

\

\paragraph{\bf Proof:} We begin as in Fefferman-Phong argument (see \cite{Fef} and also \cite{Sh1}) 
we have 
$$
\int_Q |\nabla u|^p \ge \frac C {R^{n+p}} \int_Q\int_Q |u(x)-u(y)|^p\, dxdy
$$
and 
$$
\int_Q w|u|^p = \frac {1}{R^n} \int_Q\int_Q w(x)|u(x)|^p\, dxdy.
$$
Hence,
\begin{equation}
\label{eq:FP}
\int_Q |\nabla u|^p + w|u|^p \ge  \av_Q\big[ \min (CR^{-p}, w)\big]\int_Q |u(y)|^p\, dy.
\end{equation}
Now, we use that $w \in A_\infty$.  There exists $\varepsilon>0$, independent of $Q$,  such that $E=\{x\in Q\, ; \, w(x) >\ep \av_Q w\}$ satisfies  $|E| \ge \frac 1 2 |Q|$. Hence 
$$
 \av_Q\big[ \min (CR^{-p}, w)\big]  \ge \frac 1 2 \min (CR^{-p}, \ep \av_Q w). $$
 This proves  the desired inequality when $R^p \av_Q w \le 1$. 
 
 Assume now that  $R^p\av_Qw \ge 1$. Subdivide $Q$ in a dyadic manner and stop the first time that $R(Q')^p \av_{Q'}w<1$. One obtains a collection $\{Q_i\}$ of strict dyadic subcubes of $Q$ which are maximal for the property $R_i^p\av_{Q_i} w <1$. Furthermore, since $w(x)\, dx$ is a doubling measure and as the ancestor $\widehat{Q_{i}}$ of $Q_{i}$ satisfies  $(2R_{i})^p \av_{\widehat Q_{i}}w\ge 1$, there exists $A>0$ such that 
 $R_i^p\av_{Q_i} w \ge A$. The last observation is that  the $Q_{i}$ form a disjoint covering of  $Q$ up to a set of null measure. Indeed, for almost all $x\in Q$, $\av_{Q'} w$ converges to $w(x)$, and therefore $R(Q')^2 \av_{Q'}w$ to 0, whenever $Q'$ describes the sequence of dyadic subcubes of $Q$ that contain $x$ and $R(Q')$ tends to 0. 
Hence, 
\begin{align*}
\int_Q |\nabla u|^p + w|u|^p & = \sum_i  \int_{Q_i} |\nabla u|^p + w|u|^p 
\\
& \ge C'  \sum_i  \min (R_i^{-p}, \av_{Q_i} w) \int_{Q_i} |u|^p 
\\
& \ge A C' \sum_i R_i^{-p} \int_{Q_i} |u|^p
\\
& \ge AC' \min_i \left( \frac R {R_i} \right)^p R^{-p} \int_Q|u|^p.
\end{align*}
It remains to estimate $ \min_i \left( \frac R {R_i} \right)^p$ from below.
Let $1\le \alpha <\infty$ be such that $w\in A_\alpha$. Then, 
for any cube $Q$ and measurable subset $E$ of $Q$, we have
$$
\left(\frac  {\av_{E}w} {\av_Q w} \right) \ge C   \left (\frac {|E|}  {|Q|} \right)^{\alpha-1}.
$$
Applying this to $E=Q_{i}$ and $Q$, we obtain, 
$$
\left( \frac R {R_i} \right)^p = \left(\frac {R^p \av_Q w}{R_i^p  \av_{Q_i} w }\right) \left(  \frac  {\av_{Q_i}w} {\av_Q w} \right)
\ge R^p \av_Q w \left(\frac  {\av_{Q_i}w} {\av_Q w} \right) \ge C R^p \av_Q w  \left (\frac  {R_i} R\right)^{n(\alpha-1)}.
$$
This yields $ \min_i \big( \frac R {R_i} \big)^p \ge C (R^p \av_Q w )^\beta$ with $\beta = \frac p {p+n(\alpha -1)}$ and the lemma is proved.

\begin{rem}  If $w \in B_q$ for $q>n/p$ $($as in \cite{Sh1} with $p=2)$, then $ R/R_i$ is also bounded by $C (R^p\av_Q w)^{p-n/q}$, that is $R/R_i$ is logarithmically comparable to $R^p \av_Q w$. No such thing is true if $q<n/p$.
For example, if $w(x) = |x|^{-\alpha}$ with $p<\alpha<n$ $($hence $w\in B_q$ for $q<n/\alpha<n/p)$ then  it is easy to show that $\max R/R_i$ can be unbounded. Furthermore, for all $x$ then $R^p\av_{Q(x,R)}w$ tends to 0 as $R\to +\infty$, which is not the case when $0<\alpha<p$. The case $\alpha=p$ is different in the sense that $R^p\av_{Q(x,R)}w$ tends to a non zero constant as $R\to +\infty$.
\end{rem}

\section{Definitions of the Schr\"odinger operator}

Recall that $V$ is a nonnegative locally integrable function on $\RR^n$. The definition of the Schr\"odinger operator associated to  $-\Delta+V$ is via the quadratic form method. Let
$$
\V= \{f\in L^2(\RR^n)\, ; \, \nabla f \ \& \ V^{1/2}f \in L^2(\RR^n)\}.
$$
Equipped with the norm
$$
\|f\|_{\V}= \big(\|f\|_{2}^2+ \|\nabla f\|_{2}^2 + \|V^{1/2} f\|_{2}^2\big)^{1/2}
$$it is a Hilbert space and it is known that 
$\coinf(\RR^n)$ is dense in $\V$ (\cite{Da2}).
The sesquilinear form
$$
Q(u,v) = \int_{\RR^n} \nabla u \cdot \overline{\nabla v } + V u \, \overline v 
$$
on $\V\times \V$ is bounded below and non-negative and, therefore, there exists a unique  positive self-adjoint operator
$H$ such that 
$$
\langle Hu, v\rangle = Q(u,v) \quad \forall\, u \in \D(H) \ \forall\, v\in \V.
$$
An important feature is   the pointwise domination of the resolvent by that of the Laplacian (see \cite{Da2}). This allows to define $(H+\ep)^{-1}$ on $L^p$ for $1\le p \le \infty$ and $\ep>0$ and  for any $f \in L^p$, $f\ge 0$
$$0 \le(H+\ep)^{-1} f \le (-\Delta + \ep)^{-1} f $$  

Since $\D(H)$ is dense in $\V$, there is a natural extension $\ext$ of $H$ as a bounded operator
from $\V$ to $\V'$ (not identified with $\V$). Further, for any $\ep>0$, $\ext+\ep$ is invertible but this ceases at $\ep=0$ so 
it is useful to introduce an ``homogeneous'' version of $H$ as follows: Let $\dot \V$ be the closure of 
$\coinf(\RR^n)$ under the semi-norm 
$$
\|f\|_{\dot \V}= \big( \|\nabla f\|_{2}^2 + \|V^{1/2} f\|_{2}^2\big)^{1/2}.
$$
By \eqref{eq:FP},  there is a continuous inclusion $\dot \V \subset L^2_{loc}(\RR^n)$  if  $V$ is not identically 0, which is assumed from now on,  hence, this is a norm.    
The form $Q$ is the inner product on $\dot \V$ associated to this norm so that $\dot\V$ is a Hilbert space. But if we choose not to identify $\dot \V$ and its dual, then there is 
a unique bounded and invertible operator $\schdot\colon \dot \V \to \dot \V'$ such that 
for all $u, v \in \dot\V$, $\langle \schdot u , v\rangle = Q(u,v)$. Here, $\langle\ , \ \rangle$ is the duality (sesquilinear)  form between $\dot \V'$ and $\dot \V$. Note that since $\coinf(\RR^n)$ is densely contained in $\dot \V$, this coincides with the usual duality between distributions and test functions when $v\in \coinf(\RR^n)$. By abuse, we do not distinguish the two notations, which we write as  an integral when the integrand is integrable.

In concrete terms, if $f\in \dot \V'$ there exists a unique $u\in \dot \V$ such that 
\begin{equation}
\label{eq:Htilde}
\int_{\RR^n} \nabla u \cdot {\nabla \overline v } + V u\,  \overline v = \langle f, v \rangle \quad \forall\, v\in \coinf(\RR^n).
\end{equation}
In particular, $-\Delta u + Vu = f$ holds in the sense of distributions.
There is a classical approximation procedure to obtain $u$ for nice $f$. 

\begin{lem}\label{lem:approx}
Assume that $f\in   \dot \V' \cap L^2(\RR^n)$. For $\ep>0$, let $u_{\ep}= (H+\ep)^{-1}f \in \D(H)$. Then 
$(u_{\ep})$ is a bounded sequence in $\dot \V$ which
converges strongly to $\schdot^{-1} f$.

\end{lem}

\paragraph{\bf Proof} 
By definition, 
$$
\int_{\RR^n} \nabla u_{\ep} \cdot {\nabla\overline v } + (V +\ep) u_{\ep}\, \overline v = \int_{\RR^n} f \,  \overline v \quad \forall\, v\in \V
$$
and in particular
$$
\int_{\RR^n}| \nabla u_{\ep}|^2 + (V +\ep) |u_{\ep} |^2 = \int_{\RR^n} f \, \overline {u_{\ep}}.
$$
The boundedness of $(u_{\ep})$ in $\dot\V$ follows readily using that $\vert \int_{\RR^n} f  \, \overline {u_{\ep}}\vert \le \|f\|_{\dot\V'}\|u_{\ep}\|_{\dot \V}$ and $f \in \dot \V'$. 

Let us see first the weak convergence. Let $u\in \dot \V$ be a weak limit of a subsequence $(u_{\ep_{j}})$. One can take limits in the first equation when $v\in \coinf(\RR^n)$ and we see that $u$ satisfies \eqref{eq:Htilde}. By uniqueness, $u=\schdot^{-1} f$ and $(u_{\ep})$ converges weakly to $u$. Since $f\in \dot \V'$, we have
$$
\int_{\RR^n}| \nabla u|^2 + V  |u|^2 =  {\langle f,u  \rangle }.$$
Weak convergence implies
\begin{align*}
\int_{\RR^n}| \nabla u|^2 + V  |u|^2 & \le  \liminf\int_{\RR^n}| \nabla u_{\ep}| + V  |u_{\ep} |^2
\le \limsup \int_{\RR^n}| \nabla u_{\ep}|^2 + V  |u_{\ep} |^2 
\\
&\le 
 \limsup \int_{\RR^n}| \nabla u_{\ep}|^2 + (V +\ep) |u_{\ep} |^2 
 = \limsup \int_{\RR^n} f \, \overline{u_{\ep}} = {\langle f,u  \rangle }.
 \end{align*}
Thus $\|u_{\ep}\|_{\dot \V}\to \|u\|_{\dot \V}$ and together with weak convergence, this yields strong convergence.

\begin{rem}
The continuity of the inclusion $\dot \V \subset L^2_{loc}(\RR^n)$ has two further consequences: first, we have that $L^2_{comp}(\RR^n)$, the space of compactly supported $L^2$ functions on $\RR^n$, is continuously contained in 
$\dot\V' \cap L^{2}(\RR^n)$.  Second, $(u_{\ep})$ has a subsequence converging to $u$ almost everywhere.
\end{rem}

We continue with square roots. As $H$ is self-adjoint, it has a unique square root $H^{1/2}$, which is self-adjoint with domain $\V$ and for all $u\in C_{0}^\infty(\RR^n)$, $\|H^{1/2}u\|_{2}^2= \|\nabla u\|_{2}^2+\|V^{1/2} u\|_{2}^2$.  This allows us to extend 
$H^{1/2}$ from  $\dot\V$ into $L^2$. If $S$ denotes this extension, then we have $S^*S=\schdot$ where $S^*\colon L^2 \to \dot\V'$ is the adjoint of 
$S$.

Our results are all about $\schdot$ or alternately about $H+\ep$ with a uniform control of constants with respect to $\ep>0$. By abuse, we write $\sch$ for $\schdot$ and $\sch^{1/2}$ for its extension $S$ or its adjoint $S^*$. The context will make clear which object is the right one.

\section{An $L^1$ maximal inequality}

The following result is essentially a consequence of a result of Gallou\"et and Morel \cite{GM} in the semi-linear setting or can be seen from \cite{K2}.  We present a simple proof in this situation.  We assume that $V$ is not identically 0

\begin{lem}\label{lem:l1} Let  $f \in L^{\infty}_{comp}(\RR^n)$, $f\ge 0$ and $u=\sch^{-1}f$. Then $u\ge 0$, 
 $$
\int_{\RR^n} Vu \le \int_{\RR^n} f, 
$$
 $$
\int_{\RR^n}| \Delta u| \le 2 \int_{\RR^n}f. 
$$
Furthermore, $u\in W^{1,1}_{loc}(\RR^n)$ and 
for any measurable set $E$ with bounded measure,
 $$
\int_{E}| \nabla u| \le C(n) |E|^{1/n} \int_{\RR^n}f, 
$$
and  for all compact set $K$ in $\RR^n$, 
$$
\int_{K}|  u| \le C(K,n,V) \int_{\RR^n} f.
$$
\end{lem}

\begin{rem}
In fact, more is true. If $n=1$, the estimate on $u'$ tells us that $u'$ is bounded. If $n\ge 2$, 
then $u\in W_{loc}^{1,q}(\RR^n)$ for $1\le q < \frac n{n-1}$.
\end{rem}

\paragraph{\bf Proof} For $N\ge \ep>0$, set $V_{\ep,N}= \inf( V+\ep, N)$. Let 
$f \in L^{\infty}_{comp}(\RR^n)$, $f\ge 0$ and set $u=\sch^{-1}f$, $u_{\ep}=(H +{\ep})^{-1}f$ and $u_{\ep, N}= H_{\ep,N} ^{-1} f$ where  $H_{\ep, N}$ is associated  to the potential $V_{\ep, N}$. By Lemma \ref{lem:approx} we know that $u\in L^1_{loc}(\RR^n)$ (with norm controlled by $\|f\|_{\infty}$ which is not enough). We remark that by comparison theorems   $u, u_{\ep}, u_{\ep, N}\ge 0$. Further,   $V_{\ep, N}\le V+\ep$ implies  $ u_{\ep}\le u_{\ep, N}$ and $V\le V+\ep$ implies $(u_{\ep})_{\ep>0}$ is non increasing  with $u_{\ep} \le u$.  In addition, it follows from the remark after Lemma \ref{lem:approx} that  
$u_{\ep}$ converges almost everywhere to $u$ as $\ep \to 0$.  Indeed, a subsequence already converges almost everywhere to $u$, hence the family itself by monotonicity. As a consequence, $(u_{\ep})$ converges to $u$ in $L^1_{loc}(\RR^n)$ by the monotone convergence theorem. 

Let us see the estimate for $Vu$. Since $\ep \le V_{\ep, N}\le N$, the operator $H_{\ep,N} ^{-1}$ has an integral kernel bounded by the one of $(-\Delta+\ep)^{-1}$, 
which implies that it extends to a bounded operator on  $L^p(\RR^n)$ for $1\le p \le \infty$. 
As $V_{\ep,N}$ is bounded,  $V_{\ep, N} H_{\ep,N} ^{-1}$ is a bounded operator  on $L^1(\RR^n)$ and also 
$ \Delta H_{\ep,N} ^{-1}$ by difference. As $u_{\ep, N}\in L^1(\RR^n)$ and $\Delta u_{\ep, N} \in L^{1}(\RR^n)$, an easy argument via  the Fourier transform implies that $\int_{\RR^n} \Delta u_{\ep, N} = 0$,  and so
$$
\int_{\RR^n} V_{\ep, N} u_{\ep ,N} =\int_{\RR^n} f.
$$
Next, we have $0\le V_{\ep, N} u_{\ep}\le V_{\ep, N} u_{\ep ,N}$, so $V_{\ep, N} u_{\ep}$ is integrable and by monotone convergence as $N\to \infty$,  $(V +{\ep}) u_{\ep}$ is integrable and
$$
\int_{\RR^n} (V+{\ep}) u_{\ep} \le \int_{\RR^n} f.
$$
Finally, we have
$$
\int_{\RR^n} Vu_{\ep} \le
\int_{\RR^n} (V +{\ep}) u_{\ep}  \le \int_{\RR^n} f
$$
and monotone convergence as $\ep \to 0$ yields 
$$
\int_{\RR^n} Vu \le \int_{\RR^n} f.
$$

We turn to the term with $\Delta u$. As  
$(u_{\ep})$ converges  to 
$u$  in $L^1_{loc}(\RR^n)$, $\Delta u$  is the limit of $\Delta u_{\ep}$ in $\D'(\RR^n)$.  If $h\in \coinf(\RR^n)$, then
$$
-\int_{\RR^n}  \Delta u_{\ep}\,  h =   \int_{\RR^n}  \nabla u_{\ep}\cdot \nabla h = \int_{\RR^n}  f h - \int_{\RR^n}  V u_{\ep} h - \int_{\RR^n} \ep u_{\ep} h.
$$
As $h$ has compact support, the last integral converges to $0$, hence $-\Delta u$ is equal to $ f - Vu
\in \D'(\RR^n)$ and  its $L^1$ control follows.

We turn to the gradient estimate. As  $u_{\ep}\in L^1(\RR^n)$ and $\Delta u_{\ep} \in L^1(\RR^n)$, it can be shown (see \cite[Appendix]{BBC}) that   if $E$ is a measurable subset of $\RR^n$ with bounded measure and $1\le q<\frac n{n-1}$,
$$
\int_{E}  |\nabla u_{\ep}|^q  \le C(n,q) |E|^{1- \frac {(n-1)q}{n}}     \|\Delta u_{\ep}\|_{1} ^q
$$
hence
\begin{equation}
\label{eq:graduep}
\int_{E}  |\nabla u_{\ep}|^q  \le  C(n,q) |E|^{1- \frac {(n-1)q}{n}}   \|f\|_{1}^q.
\end{equation}

Next,  recall that if $Q$ is a cube, by \eqref{eq:FP} 
we have, 
$$
 \av_Q\big[ \min (CR^{-1}, V)\big]\int_Q  u_{\ep}  \le \int_Q |\nabla u_{\ep}| + Vu_{\ep} $$
 hence
\begin{equation}
\label{eq:uep}
\int_Q  u_{\ep}   \le C(Q,n, V) \|f\|_{1}.
\end{equation}
It follows easily from these two estimates and Poincar\'e inequality  that $u_{\ep}\in W^{1,q}_{loc}(\RR^n)$ for $1\le q<\frac n{n-1}$ and is bounded in that space. 
Thus, for any 
$1<q<\frac n{n-1}$,  $u \in W^{1,q}_{loc}(\RR^n)$ by taking weak limits.  The estimate \eqref{eq:graduep} passes to the liminf and by H\"older becomes true for $q=1$. The estimate \eqref{eq:uep} also passes to the limit by convergence in $L^1_{loc}(\RR^n)$. This finishes the proof.

 \

Let $B=\{ u\in L^1_{loc}(\RR^n)\, ; \, \Delta u \in L^1(\RR^n),  Vu\in L^1(\RR^n)\}$ equipped with the topology defined by the semi-norms for $L^1_{loc}(\RR^n)$,  $\|\Delta u \|_{1}$ and $\|Vu\|_{1}$. We have obtained

\begin{theo}\label{th:L1}  The operator $\sch^{-1}$ \textit{a priori} defined on  
$L^2_{comp}(\RR^n) $ extends to  a bounded operator  from $L^1(\RR^n)$ into $B$. Denoting again $\sch^{-1}$ this extension, $V\sch^{-1}$ is a  positivity-preserving contraction on $L^{1}(\RR^n)$ and $\frac 1 2 \Delta \sch^{-1}$ a contraction on $L^{1}(\RR^n)$.\end{theo}

\begin{pro}\label{pro:uniqueness} Let $f\in L^1(\RR^n)$.
There is uniqueness of solutions for the  equation $-\Delta u + Vu =f $ in the class $L^1(\RR^n)\cap B$. In particular, if $u \in C_{0}^\infty(\RR^n)$ and $f=-\Delta u + Vu$, then 
$u=\sch^{-1}f$. 
\end{pro}

\paragraph{\bf Proof} Since $-\Delta u + Vu =0$, then for $\ep>0$  we have $-\Delta u + Vu + \ep u = \ep u $. As $u\in L^1(\RR^n)$, we can write 
$|u| \le (-\Delta +\ep)^{-1}(\ep | u|)= (-\ep^{-1} \Delta + 1)^{-1}| u|$. Taking limits as $\ep \to 0$ proves that $u=0$.

\begin{rem}
We have obtained existence in $B$ and  uniqueness  in $B\cap L^1(\RR^n)$. Following \cite{GM}, one can show  uniqueness  in a larger space if $n\ge 3$ and under some conditions on $V$ if $n\le 2$. We do not need such refinements here.  
\end{rem}

\begin{cor}\label{cor:L1}
Equation \eqref{eq:L1} holds.
\end{cor}

\paragraph{\bf Proof}    If $u \in C_{0}^\infty(\RR^n)$ and $f=-\Delta u + Vu$, then 
$Vu=V\sch^{-1}f$ and $\Delta u = \Delta  \sch^{-1} f$ by the proposition above.  Applying Theorem \ref{th:L1} proves that 
$\|Vu \|_{1} \le  \|-\Delta u + Vu \|_{1}$ and $\|\Delta u\|_{1}\le 2 \|-\Delta u + Vu \|_{1}$.

\section{$L^p$ maximal inequalities}\label{sec:lp}

The main sledge hammer  is the following  criterion for $L^p$ boundedness (\cite{AM1}).  A slightly weaker version appears in Shen \cite{Sh2}. 
 
 \begin{theo} \label{theor:shen}
Let $1\le p_{0} <q_0\le \infty$. Suppose that $T$ is a bounded
sublinear operator  on $L^{p_{0}}(\RR^n)$. Assume that there exist
constants $\alpha_{2}>\alpha_{1}>1$, $C>0$ such that
\begin{equation}\label{T:shen}
\big(\av_{Q} |Tf|^{q_0}\big)^{\frac1{q_0}}
\le
C\, \bigg\{ \big(\av_{\alpha_{1}\, Q}
|Tf|^{p_0}\big)^{\frac1{p_0}} +
(S|f|)(x)\bigg\},
\end{equation}
for all cube $Q$, $x\in Q$ and  all $f\in L^{\infty}_{comp}(\RR^n)$ with 
support in $\RR^n \setminus \alpha_{2}\, Q$, where $S$ is a positive operator.
Let $p_{0}<p<q_{0}$. If $S$ is bounded on $L^p(\RR^n)$, then, there is a constant $C$ such that
$$
\|T f\|_{p}
\le
C\, \|f\|_{p}
$$
for all    $f\in L_{comp}^{\infty}(\RR^n)$.\end{theo}

Note that in this statement, $f$ can be valued in a Banach space and $|f|$ denotes its norm.
Also  the space $L_{comp}^{\infty}(\RR^n)$ can be replaced by $C_{0}^\infty(\RR^n)$.

\
 
Fix an open set $\Omega$. By a weak solution of $-\Delta u + Vu=0$ in $\Omega$, we mean $u\in L^1_{loc}(\Omega)$ with       $V^{1/2}u,   \nabla u \in L^2_{loc}(\Omega)$ and the equation holds in the distribution sense on $\Omega$. Remark that by Poincar\'e's inequality if $u$ is a weak solution, then $u \in L^2_{loc}(\Omega)$. A subharmonic function on $\Omega$ is a function $v\in L^{1}_{loc}(\Omega)$ such that $\Delta v \ge 0$ in $\D'(\Omega)$. It should be observed that if $u$ is a weak solution in $\Omega$ of 
$-\Delta u + Vu=0$  then 
\begin{equation}
\label{eq:sub}
 \Delta |u|^2 = 2 V|u|^2 + 2 |\nabla u |^2
\end{equation}
 and in particular, $|u|^2$ is  a nonnegative  subharmonic  function in $\Omega$. 

The main technical lemma is interesting on its own right. It states that a form of the  mean value inequality for subharmonic functions still holds if the Lebesgue measure is replaced 
by a weighted measure of Muckenhoupt type. More precisely, 

\begin{lem}\label{lem:subw}
Assume $w\in A_{\infty}$ and $f$ is a nonnegative  subharmonic function in  $\Omega$, $Q$ is a cube in $\RR^n$ with $\overline{2Q} \subset \Omega$,  $1 < \mu \le 2$ and 
 $0<s<\infty$. Then for some $C$ depending on the $A_\infty$ constant of $w$, $s$, $\mu$ (and independent of $f$ and $Q$) and  almost all $x\in Q$,  we have
$$f(x) \le \bigg(\frac {C}{w(\mu Q)} \int_{\mu Q} w f^s\bigg)^{1/s}.$$
\end{lem}

Here $w(E)=\int_{E}w$. 
  As $A_{\infty}$ weights have the doubling property we have $\av_{\mu Q}w\sim \av_{Q}w$ and the inequality above
rewrites (the notation $\sup$ meaning essential supremum)
\begin{equation}
\label{eq:subbis}
\big(\av_{Q}w\big) \big(\sup_{Q}f^s\big) \le C \av_{\mu Q}( w f^s).
\end{equation}

\paragraph{\bf Proof} This is a consequence of a result of S. Buckley \cite{Bu}. We give the proof
for the convenience of the reader. Since $w\in A_{\infty}$, there is $t<\infty$ such that  $w\in A_{t}$. Hence
for any nonnegative measurable function $g$, we have 
$$
\av_{\mu Q} g \le C \Big(\frac1{w(\mu Q)} \int_{\mu Q} w  g^t \Big)^{1/t} = C \big( \av_{\mu Q} (w g^t )\big)^{1/t} \big( \av_{\mu Q}  w\big)^{-1/t} .  
$$
But  subharmonicity of $f$ in $\Omega$ implies for almost all $x\in Q$  and all $0<r<\infty$
\begin{equation}
\label{eq:subharmonic}
f(x) \le C_{r, \mu}\, \big(\av_{\mu Q} f^r\big)^{1/r}
\end{equation}
 (See \cite{FS}. It can also be obtained from classical facts on weak reverse H\"older weights \cite{IN}). Applying this  with $r=s/t$   yields
$$
f (x)\le C \big(\av_{\mu Q} f^{s/t}\big)^{t/s} \le  C  \Big( \av_{\mu Q} (w f^s )\Big)^{1/s} \big( \av_{ \mu Q}  w\big)^{-1/s} .
$$

\begin{cor}\label{lem:rh2q} 
Let  $w\in B_{r}$ for some $1<r\le \infty$ and let  $0<s<\infty$.  Then there is  $C\ge 0$  depending only on the $B_r$ constant of $w$, $s$, $\mu$  such that  for any  cube $Q$  and any nonnegative  subharmonic function  $f$ in  a neighborhood of $\overline{2Q}$  we have  for all $1< \mu \le 2$
$$
\big(\av_{Q} (w f^s)^{r} \big)^{1/r} \le C \, \av_{\mu Q} (w f^s).
$$
\end{cor}

\paragraph{\bf Proof}  We have
$$
\big(\av_{Q} (w f^s)^r\big)^{1/r} \le  \big(\av_{Q} w^r \big)^{1/r} \sup_{Q} f^s \le 
C\, \big(\av_{Q}w\big) \sup_{Q} f^s \le C \, \av_{\mu Q} (w f^s).
$$
The second inequality uses the $B_{r}$ condition on $w$ and the last inequality is \eqref{eq:subbis}.

\

Let us come back to the Schr\"odinger operator.

\begin{theo} Let $V\in B_{q}$ for some $1<q\le \infty$. Then there is $r>q$ (or $r=\infty$ if $q=\infty$) such that
$V\sch^{-1}$ and $\Delta\sch^{-1}$, defined on $L^1(\RR^n)$ from Theorem \ref{th:L1}, extend to  bounded operators on $L^p(\RR^n)$ for
$1<p<r$.

\end{theo}

\paragraph{\bf Proof} By difference, it suffices to prove the theorem for $V\sch^{-1}$. We know that this is a bounded operator on $L^1(\RR^n)$. Let $r>q$ be given by  self-improvement of the reverse H\"older inequalities of $V$. Fix a cube $Q$ and let $f\in L^\infty(\RR^n)$ with compact support contained in $\RR^n \setminus 4Q$. Then $u=\sch^{-1} f$ is well-defined in    $\dot\V$ and is a weak solution of $-\Delta u + Vu=0$ in $4Q$. Since $|u|^2$ is subharmonic, the above corollary  applies  with $w=V$, $f=|u|^2$ and $s=1/2$.   It  yields \eqref{T:shen} with $T=V\sch^{-1}$,  $p_{0}=1$, 
$q_{0}=r$, $S=0$, $\alpha_{1}=2$ and $\alpha_{2}=4$. Hence, $T$ is bounded on $L^p(\RR^n)$ for $1<p<r$ by Theorem \ref{theor:shen}. 

\

\paragraph{\bf Proof of Theorem \ref{th:1} } 
   Let $u\in \coinf(\RR^n)$ and $f=-\Delta u + Vu$. We know that $u=\sch^{-1}f$ by Proposition \ref{pro:uniqueness}. Now, using the hypothesis  $V\in B_{q}$, we have  bounded extensions on $L^p(\RR^n)$  
of $V\sch^{-1}$ and $\Delta \sch^{-1}$ for $1<p<q+\ep$ for some $\ep>0$ depending on $V$. 
We conclude that  $\|Vu \|_{p} +\|\Delta u\|_{p}\lesssim \|f \|_{p}.$

\section{Complex interpolation}\label{sec:interpolation}

We shall use complex interpolation to obtain   item 1  of Theorem \ref{th:2}, relying on 
the following result due to Hebisch \cite{H}.

\begin{pro}
Let $V$ be a nonnegative locally integrable function on $\RR^n$. Then for all $y\in \RR$,  $H^{iy}$  has a bounded extension on $L^p(\RR^n)$, $1<p<\infty$, and for fixed $p$ its operator  norm does not exceed
$C(\delta,p )e^{\delta |y|}$ for all $\delta >0$.  
\end{pro}

Here, $H^{iy}$ is defined  as  a bounded operator on $L^2(\RR^n)$ by functional calculus. For $V=0$, this is  standard result for the singular integral operator $(-\Delta)^{iy}$. Actually, the operator norm can be improved but we do not need sharp estimates.

\begin{lem}
The space $\D=\mathcal{R}(H) \cap L^1(\RR^n)\cap L^\infty(\RR^n)$ is dense in $L^p(\RR^n)$ for $1<p<\infty$. 
\end{lem}

\paragraph{\bf Proof} It suffices to show that $\mathcal{R}(H) \cap L^1(\RR^n)\cap L^\infty(\RR^n)$ is dense in $L^1(\RR^n)\cap L^\infty(\RR^n)$ for the $L^p(\RR^n)$ norm. Let $f\in L^1(\RR^n)\cap L^\infty(\RR^n)$. Since $f\in L^2(\RR^n)$,  for $\ep>0$, $f_{\ep}=
H(H+\ep)^{-1}f  \in \mathcal{R}(H)$. Also $f_{\ep}= f - \ep(H+\ep)^{-1}f$. Thus $f_{\ep} \in L^1(\RR^n)\cap L^\infty(\RR^n)$  as $|(H+\ep)^{-1}f| \le (-\Delta+\ep)^{-1}|f|$
and the kernel of $(-\Delta+\ep)^{-1}$ is integrable. It remains to see that $f_{\ep}$ converges to $f$ in $L^p(\RR^n)$. But again $|f-f_{\ep}| \le \ep (-\Delta+\ep)^{-1}|f|$ and the latter expression is easily seen to converge to $0$ in $L^p$ as $\ep$ tends to 0.

\

We now prove the boundedness of $\nabla H^{-1/2}$ and $V^{1/2} H^{-1/2}$ on $L^p(\RR^n)$ for $1<p<2(q+\ep)$, which is half of  item 3 of Theorem \ref{th:2}.
Let $f\in \D$, $g\in C_{0}^\infty(\RR^n)$.  We define for $z \in S=\{ x+iy\in \CC\, ;\, 0 \le x\le 1, y\in \RR\}$,
$$
A(z)= \langle (-\Delta)^z H^{-z} f, g\rangle 
$$
We shall use the Stein interpolation theorem for families of operators (see \cite{SW}).

Observe that for all $z\in S$,  $(-\Delta)^{\bar z} g \in L^2$ with $\|(-\Delta)^{\bar z} g\|_{2}\le C \|g\|_{H^2}$ (the Sobolev space of order 2). Since $f\in \mathcal{R}(H)$, $f=H\tilde f$ with $M=\|\tilde f\|_{2}+ \|H\tilde f\|_{2}<\infty$. Hence, 
$$
\|H^{-z}f\|_{2}= \|H^{1-z} \tilde f\|_{2} \le \|H^{-iy}\|_{2,2} \|H^{1-x}\tilde f\|_{2} \le C({\delta }) e^{\delta |y|}  M.
$$
Thus $|A(z)| \le  C_{\delta } e^{\delta |y|}  M \|g\|_{H^2}$. It follows that $A$ satisfies the admissible growth condition. It is not difficult to establish continuity on $S$ and analyticity on ${\rm Int}\, S$ of $A$.  Then, for $z=iy$ and $1<p<\infty$, we have
$$
|A(iy)| \le  \|H^{-iy}f\|_{p} \|(-\Delta)^{-iy} g\|_{p'} \le C({\delta, {p}}) e^{\delta |y|} \|f\|_{p}\|g\|_{p'}.
$$
And for $z=1+iy$ and $1<p<q+\ep$, 
$$
|A(1+ iy)| \le  \|\Delta H^{-1} H^{-iy}f\|_{p}\|(-\Delta)^{-iy} g\|_{p'}  \le \|\Delta H^{-1}\|_{p,p} 
C({\delta, {p}}) e^{\delta |y|} \|f\|_{p}\|g\|_{p'}.
$$
Thus, for $z=1/2$ and $1<p<2(q+\ep)$, we obtain
$$
|A(1/2)| \le  C(p) \|f\|_{p}\|g\|_{p'}.
$$
We conclude by a density argument that $(-\Delta)^{1/2} H^{-1/2}$ is bounded on $L^p(\RR^n)$ for $1<p<2(q+\ep)$.

Similarly, for $f\in \D$, $g\in C_{0}^\infty(\RR^n)$, we define for $z \in S=\{ x+iy\in \CC\, ;\, 0 \le x\le 1\}$
and fixed $N>0$, 
$$
B(z)= \langle V_{N}^z H^{-z} f, g\rangle, 
$$
with $V_{N} =\inf  (V, N)$. 
Then, 
$$
|B(z)| \le  C({\delta }, q') e^{\delta |y|} M (\|g\|_{2} + \|V_{N} g\|_{2}), 
$$
hence $B$ has the admissible growth condition. It is also clearly continuous on $S$ and analytic on ${\rm Int}\, S$. Then, for $z=iy$ and $1<p<\infty$. 
$$
|B(iy)| \le  \|H^{-iy}f\|_{p} \|V_{N}^{-iy} g\|_{p'} \le C({\delta, {p}}) e^{\delta |y|}\|f\|_{p}\|g\|_{p'}.
$$
And for $z=1+iy$ and $1<p<q+\ep$, 
$$
|A(1+ iy)| \le  \| V_{N} H^{-1} H^{-iy}f\|_{p}\|V_{N}^{-iy} g\|_{p'}  \le \|V H^{-1}\|_{p,p} 
C({\delta, {p}}) e^{\delta |y|} \|f\|_{p}\|g\|_{p'}.
$$
where we used that $\|V_{N} H^{-1}\|_{p,p} \le \|V H^{-1}\|_{p,p}$ as $0\le V_{N}\le V$ almost everywhere.
Thus, for $z=1/2$ and $1<p<2(q+\ep)$, we obtain
$$
|B(1/2)| \le  C(p) \|f\|_{p}\|g\|_{p'}.
$$
We conclude by a density argument that $V_{N}^{1/2} H^{-1/2}$ is bounded on $L^p(\RR^n)$ for $1<p<2(q+\ep)$ with a bound that is uniform with respect to $N$. By monotone convergence, this yields the $L^p(\RR^n)$ boundedness  of $V^{1/2} H^{-1/2}$ in the same range. 

To finish the proof, fix $1<p<2(q+\ep)$. Let $u\in C_{0}^\infty(\RR^n)$. The only thing to establish is 
$\|\nabla u \|_{p}+\|V^{1/2}u\|_{p} \le C(p) \|H^{1/2}u\|_{p}$. Since $u \in \V$, $f=H^{1/2}u$ is well-defined. We assume that $f\in L^p(\RR^n)$,  otherwise there is nothing to prove. Then, by Calder\'on-Zygmund theory and the above, 
$$
\|\nabla u \|_{p}+\|V^{1/2}u\|_{p} \le C(p) \|(-\Delta)^{1/2} H^{-1/2} f\|_{p} + \|V^{1/2}H^{-1/2} f\|_{p} \le C'(p) \|f\|_{p}$$
and the proof is finished.

\begin{rem}
We remark that this interpolation argument gives also a proof of the $L^p$ boundedness of  $\nabla H^{-1}$ and $V^{1/2}H^{-1/2}$ for $1<p<2$ for all non zero $V\in L^1_{loc}(\RR^n)$.
\end{rem}

\section{Reverse Riesz transforms}\label{sec:item3}

 This section is concerned with the proof  of  Theorem \ref{th:2},  item 2.
We first want to show that there exists $C>0$ depending only on the $A_{\infty}$ constant of $V$ such that  for all $\alpha>0$ and $f\in C_{0}^\infty(\RR^n)$ then 
\begin{equation}
\label{eq:wt}
|\{x\in \RR^n\, ; \,  |H^{1/2} f(x)| > \alpha\}| \le \frac{C}{\alpha} \int |\nabla f| + V^{1/2}|f|.
\end{equation}
First, it is not too hard to show that if $ \ep \le V \le N $ for some 
$N> \ep>0$ then this inequality holds with $C$ depending on $\ep,N$ (in fact, the next argument gives this also). Let $C_{1}$ be the best constant in this inequality with $V$ replaced by $V_{\ep,N}=\min(V+\ep,N)$. We want to show that $C_{1}$ is bounded independently of $\ep$ and $N$.  Assume it is the case, then for $\ep,N>0$, all $\alpha>0$ and $f\in C_{0}^\infty(\RR^n)$
\begin{align*}
|\{x\in \RR^n\, ; \,  |(-\Delta + V_{\ep,N})^{1/2} f(x)| > \alpha\}| &\le \frac{C_{1}}{\alpha} \int |\nabla f| + V_{\ep,N}^{1/2}|f|
\\
& \le \frac{C_{1}}{\alpha} \int |\nabla f| + (V +\ep)^{1/2}|f|.
\end{align*} 
Now, it is easy to show that $(-\Delta+V_{\ep,N})^{1/2} f$ converges in $L^2$ to $(H+\ep)^{1/2}f$ hence up to extraction of a subsequence, the above inequality passes to the limit as  $N\to +\infty$. Then,  as $f\in C_{0}^\infty(\RR^n) \subset \V=\D(H^{1/2})$, $(H+\ep)^{1/2}f$ converges to $H^{1/2}f$  in $L^2(\RR^n)$ by functional calculus as $\ep$ tends to 0 and  we obtain \eqref{eq:wt} with  $C=C_{1}$.

 Remark that if $V\in A_{\infty}$, then
for all $N>\ep>0$,  $V_{\ep,N}$ is also in $A_{\infty}$ with constants that are uniform with respect to 
$\ep$ and $N$. So as long as we only use the $A_{\infty}$ information, we are safe. Therefore, we assume that $ \ep \le V \le N$ but we do not use this information quantitatively. 

We also define $C_{p}$ as the best constant $C$ such that for $1<p<2$ and $f\in C_{0}^\infty(\RR^n)$
$$
\|H^{1/2} f\|_p \le C_p (\|\nabla f\|_p + \|V^{1/2}f\|_p).
$$
By extension, we can take $f$ to be in the closure of $C_{0}^\infty(\RR^n)$ for the norm defined by the right hand side. Since $V$ is bounded below and above, this is the usual Sobolev space $W^{1,p}(\RR^n)$. 

We know that $C_{2}=1$. We shall prove that  for some  numbers $C, M$ under control, we have
\begin{equation}
\label{eq:wt1}
C_{1}\le C C_{2}^2 + M= C+M. 
\end{equation}
This will require the use of a specific Calder\'on-Zygmund decomposition on $f$ adapted to level sets of $|\nabla f| + V^{1/2}|f|$.

The Marcinkiewicz  interpolation theorem would give us
\begin{equation}
\label{eq:marcin}
C_{p}\lesssim C_{1}^{ \frac 2p-1}
\end{equation}
provided it applies. 
But it is not known whether the spaces defined by the seminorms $\|\nabla f\|_q + \|V^{1/2}f\|_{q}$, $1\le q\le 2$,  interpolate 
by the real method. If we use the assumption 
$\ep \le V \le N$, then we may interpolate but  the constants would depend on $\ep, N$.  Instead, we   prove  \eqref{eq:marcin}  by adapting Marcinkiewicz theorem argument using again  our Calder\'on-Zygmund decomposition.

\begin{lem}\label{lemmaCZD} Let $n\ge 1$,  $1\le p< 2$, $V\in A_{\infty}$ and $f\in C_{0}^\infty(\RR^n)$, hence $\|\nabla f\|_{p} +\|V^{1/2}f\|_{p} <\infty$. Let $\alpha>0$. Then, one can find a collection of cubes $(Q_i)$, functions $g$ and $b_i$  such that 
\begin{equation}\label{eqcsds1}
f= g+\sum_i b_i  \end{equation}
and the following properties hold:
\begin{equation}\label{eqcsds2}
\|\nabla g\|_2 + \|V^{1/2}g\|_{2} \le C\alpha^{1-p/2}(\|\nabla f\|_{p} +\|V^{1/2}f\|_{p})^{p/2}, 
\end{equation}
\begin{equation}\label{eqcsds3}
\supp\,b_i  \in {Q_i}\ \text{and} \ \int_{Q_i} |\nabla b_i|^p +  R_{i}^{-p}| b_i|^p \le C\alpha^p |Q_i|, \end{equation}
 \begin{equation}\label{eqcsds4}
\sum_i |Q_i| \le C\alpha^{-p} \int_{\RR^n} |\nabla f|^p +  | V^{1/2}f |^p, \end{equation}
\begin{equation}\label{eqcsds5}
\sum_i {\bf 1}_{Q_i} \le N, \end{equation}
where 
 $N$ depends only on dimension and $C$ on dimension,  $p$ and the $A_{\infty}$ constant of $V$. Here, $R_{i}$ denotes the sidelength of $Q_{i}$ and gradients are taken in the sense of distributions in $\RR^n$.
\end{lem}

 We remark that the decomposition is on $f$ while the control is on $|\nabla f|^p +  | V^{1/2}f |^p$. 
 
 \
 
\paragraph{\bf \bf Proof} Let $\Omega$ be the open set $ \{x \in \RR^n; M(|\nabla f|^p+ | V^{1/2}f |^p)(x) >\alpha^p\}$ where $M$ is the uncentered maximal operator over cubes of $\RR^n$.
If $\Omega$ is empty, then set $g=f$ and $b_{i}=0$. Otherwise, the maximal theorem gives us
\begin{equation*} 
|\Omega| \le C\alpha^{-p} \int_{\RR^n} |\nabla f|^p+| V^{1/2}f |^p.  \end{equation*}

Let $(Q_i)$ be a Whitney decomposition of $\Omega$ by dyadic cubes:  $\Omega$ is the disjoint union of the $Q_i$'s, 
the cubes $2Q_i$ are contained in  $\Omega$ and have the bounded overlap property, but the cubes $4Q_i$ intersect $F=\RR^n\setminus \Omega$.\footnote{In fact, the factor 2 should be some $c=c(n)>1$ explicitely given in \cite{St1}. We use this convention to avoid too many irrelevant constants.}
As usual, $\lambda Q$ is the cube co-centered with $Q$ with sidelength $\lambda$ times that of $Q$. Hence \eqref{eqcsds4} and 
\eqref{eqcsds5} are satisfied by the cubes $2Q_i$. 
We remark that since $V\in A_{\infty}$, we have $V^{p/2} \in A_{\infty}$ when $1\le p\le 2$ (see Section \ref{sec:Ainfty}).  Hence we have  by Lemma \ref{lem:FP}
$$
\int_{2Q_{i}} |\nabla f|^p + |V^{1/2}f|^p \ge C \min (\av_{2Q_{i}} V^{p/2}, R_{i}^{-p} )\int_{2Q_{i}}|f|^p. 
$$ 
We declare $Q_{i}$ of type 1 if 
  $\av_{2Q_{i}} V^{p/2} \ge {R_{i}^{-p}}$
 and of type 2 if $\av_{2Q_{i}} V^{p/2} <  {R_{i}^{-p}}$.

Let us now define the functions $b_i$. Let $(\calX_i)$ be a partition of unity on $\Omega$
associated to the covering $(Q_i)$ so that for each $i$, $\calX_i$ is a $C^1$ function supported in $2Q_i$ with $\|\calX_i\|_\infty +
R_i \|\nabla \calX_i\|_\infty \le c(n)$.   
Set
$$
b_i = \begin{cases} f\calX_i, & \mathrm{if}\ Q_{i}\ \mathrm{is\ of\ type\ 1}, 
\\
 (f-\av_{2Q_{i}}f)\calX_i, & \mathrm{if}\ Q_{i}\ \mathrm{is \ of \ type\  2}.
\end{cases}
$$

If $Q_{i}$ is of type 2, then it is a direct consequence of the $L^p$-Poincar\'e inequality that 
$$
\int_{2Q_i} |\nabla b_i|^p +  R_{i}^{-p}| b_i|^p \le C \int_{{2Q_{i}}}  |\nabla f|^p.
$$
As $\int_{4Q_{i}} |\nabla f|^p \le \alpha^p  |4 Q_i|$ we get the desired inequality in  \eqref{eqcsds3}.

If $Q_{i}$ is of type 1, 
$$
\int_{2Q_i}   R_{i}^{-p}| b_i|^p \le  \int_{2Q_i}   R_{i}^{-p}| f |^p \le C\int_{2Q_{i}} |\nabla f|^p + |V^{1/2}f|^p.
$$
As the same integral but on $4Q_{i}$ is controlled by $\alpha^p|4 Q_i|$ we get 
$\int_{2Q_i}   R_{i}^{-p}| b_i|^p \le C\alpha^p| Q_i|$.  Since $\nabla b_i = \calX_i \nabla f  + f\nabla\calX_{i}$ we obtain the same bound for
 $\int_{2Q_i} |\nabla b_i|^p$. 
 
 Set $g=f-\sum b_{i}$ where the sum is over both types of cubes and is locally finite by 
 \eqref{eqcsds5}.  It is clear that
 $g=f$ on $F=\RR^n\setminus \Omega$ and $g=\sum{}^2\ (\av_{2Q_{i}} f)\ \calX_{i}$ on $\Omega$, where $\sum{}^j$ means that we are summing over cubes of type $j$.    Let us prove \eqref{eqcsds1}.
 
 First, by the differentiation  theorem, $V^{1/2}|f| \le \alpha$ almost everywhere on $F$.  Next, since
 $V \in A_{\infty}$ and $p<2$ implies $V^{p/2} \in B_{{2/p}}$ (see Section \ref{sec:Ainfty}) and 
 $\av_{2Q_{i}} V \le C ( \av_{2Q_{i}} V^{p/2})^{2/p}$. 
 Hence 
 $$
 \int_{\Omega}  V|g|^2 \le \sum{}^2 \int_{2Q_{i}} V  |\av_{{2Q_{i}}} f|^2
 \le  C \sum{}^2   \left((\av_{2Q_{i}} V^{p/2} ) |\av_{{2Q_{i}}} f|^p\right)^{2/p}  \ |Q_i|.
$$
Now, by construction of the  type 2 cubes and the $L^p$ version of Fefferman-Phong inequality,   
$$
(\av_{2Q_{i}} V^{p/2})  |\av_{{2Q_{i}}} f|^p \le C\av_{2Q_{i}}( |\nabla f|^p +  | V^{1/2}f |^p) \le C\alpha^p.
$$
 Hence, 
$$
 \int_{\Omega}  V|g|^2  \le C \sum{}^2\ \alpha^2    \ |Q_i| \le C\alpha^{2-p} \int_{\RR^n}|\nabla f|^p +  | V^{1/2}f |^p. 
 $$
 Combining the estimates on $F$ and $\Omega$, we obtain the desired bound for $ \int_{\RR^n}  V|g|^2$. 
 We finish the proof by  estimating  $ \|\nabla g\|_{\infty}$ and $\|\nabla g\|_{p}$. First, it is
 easy to see that the inequality $\|b_{i}\|_{p}^p \le C \alpha^p R_{i}^p |Q_i|$  together with the fact that Whitney cubes have sidelength comparable to their distance to the boundary, imply that 
 $\sum b_{i}$ converges in the sense of distributions in $\RR^n$ (not just  in $\Omega$, which is a trivial fact!), hence $\nabla g= \nabla f - \sum \nabla b_{i}$. 
It follows from the $L^p$ estimates on $\nabla b_{i}$ and the bounded overlap property that
$$
\left\|\sum \nabla b_{i}\right\|_{p} \le C(\|\nabla f\|_{p} + \|V^{1/2}f\|_{p}),$$
therefore the same estimate holds for $\|\nabla g\|_{p}$.
  Next, a computation of the sum $\sum \nabla b_i$ leads us to 
 $$
 \nabla g = {\bf 1}_{F} (\nabla f)   +  \sum{}^2\ (\av_{2Q_{i}} f)\ \nabla \calX_{i}.
 $$
By definition of $F$ and the differentiation theorem, $|\nabla g|$ is bounded by $\alpha$ almost everywhere on $F$.  It remains to control $\|h_{2}\|_{\infty}$ where $h_{2}= \sum{}^2\ (\av_{2Q_{i}} f)\ \nabla \calX_{i}$. Set $h_{1}=\sum{}^1 \ (\av_{2Q_{i}} f)\ \nabla \calX_{i}$. By already seen arguments  for type 1 cubes,  $|\av_{2Q_{i}} f| \le C\alpha R_{i} $. Hence,
$|h_{1}| \le C\sum{}^1\ {\bf 1}_{2Q_{i}} \alpha\le CN\alpha$ and it suffices to show that $h=h_{1}+h_{2}$ is bounded by $C\alpha$. To see this, observe that 
$\sum_i \calX_i(x)= 1$ on $\Omega$ and 0 on $F$. Since it is  a locally finite sum we have $\sum_i  \nabla\calX_i(x)=0$ for $x \in \Omega$.  Fix $x \in \Omega$. Let
$Q_j$ be the Whitney cube containing
$x$ and let $I_x$ be the set of indices $i$ such that $x \in 2Q_i$. We know that $\sharp I_x \le N$.  Also for $i \in I_x$ we have that 
$C^{-1}R_i \le  R_j \le CR_i$ (see \cite{St1}).  Therefore, we may write 
$$
|h(x)| = \left|\sum_{i \in I_x} (\av_{2Q_{i}} f- \av_{2Q_{j}} f) \nabla\calX_i(x)\right| \le C \sum_{i \in I_x} |\av_{2Q_{i}} f- \av_{2Q_{j}} f| R_i^{-1}.$$
But $2Q_{i}$ and $2Q_{j}$ are contained in $CQ_{j}$ for some $C>4$ independent of $j$. Hence, the Poincar\'e inequality and the definition of $Q_{j}$ yields
$$|\av_{2Q_{i}} f- \av_{2Q_{j}} f| \le
 CR_{j} (\av_{CQ_{j}} |\nabla f|^p)^{1/p} \le CR_{j} \alpha.
 $$ We have finished the proof.

\

\paragraph{\bf Proof of  item 3 in Theorem \ref{th:2}} First, we prove \eqref{eq:wt1}. Let $f
\in C_0^\infty(\RR^n)$.  We use the following resolution of $H^{1/2}$:
$$
H^{1/2}f= c\int_0^\infty H e^{-t^2 H} f \, dt
$$
where $c=2\pi^{-1/2}$ is forgotten from now on. It suffices to obtain the 
result for the truncated integrals $\int_\ep^R\ldots$ with bounds independent of $\ep,R$, 
and then to let
$\ep\downarrow 0$ and $R\uparrow \infty$. 
For the truncated integrals, all the calculations are justified. We thus consider that $H^{1/2}$ is one of the truncated integrals but we still write the limits as 0 and $+\infty$ to simplify the exposition. 

Apply the Cal\-de\-r\'on-Zygmund decomposition of Lemma \ref{lemmaCZD}  with $p=1$
to $f$  at height $\alpha$ and write $f=g+\sum_i b_i$.

Concerning $g$, we have
 \begin{equation*}\label{eq22}
\bigg|\bigg\{x \in \RR^n; |H^{1/2}g(x)| >\frac \alpha 3 \bigg\}\bigg| \le 
\frac{9}{\alpha^2}\int
|H^{1/2} g|^2 \le \frac{9}{\alpha^2}\int |\nabla g|^2 + V|g|^2 \le 
\frac{C}{\alpha}\int |\nabla f| +  |V^{1/2}f|
\end{equation*}
where we used   
 \eqref{eqcsds2}.
 
 The argument to estimate  $H^{1/2}b_i$ will  use the Gaussian upper bounds of the kernels of $e^{-tH}$ which are valid for all potentials $V\ge 0$. Let $r_i=2^k$ if $2^k \le R_i < 2^{k+1}$ ($R_{i}$ is the sidelength of $Q_{i}$) and
set 
$T_i= \int_0^{r_i} He^{-t^2 H}  \, dt$ and
 $U_i= \int_{r_i}^\infty He^{-t^2 H}  \, dt$. It is enough to 
estimate $$A=|\{x \in \RR^n; |\sum_i T_ib_i(x)| >\alpha/3\}|$$ and 
$$B=|\{x \in \RR^n; |\sum_i U_ib_i(x)| >\alpha/3\}|.$$

First,  
$$A \le |\cup_i 4Q_i| + \bigg|\bigg\{x \in \RR^n \setminus \cup_i  4Q_i ; 
\bigg|\sum_i T_ib_i(x)\bigg| >\frac \alpha 3\bigg\}\bigg|,$$
and by \eqref{eqcsds4}, $|\cup_i 4Q_i| \le \frac{C}{\alpha}\int |\nabla f|+  |V^{1/2}f|$.

For the other term, we have
$$
\bigg|\bigg\{x \in \RR^n \setminus \cup_i 4Q_i ; \bigg|\sum_i
T_ib_i(x)\bigg| >\frac \alpha 3\bigg\}\bigg|
\le
\frac{C}{\alpha^2}\int \bigg|\sum_i h_i\bigg|^2$$ with $h_i = {\bf
1}_{(4Q_i)^c}|T_ib_i|$. To estimate the $L^2$ norm, we dualize against
$u\in L^{2}(\RR^n)$ with $\|u\|_{2}=1$:
$$
\int| u| \sum_i h_i = \sum_i\sum_{j=2}^\infty A_{ij}
$$  
where 
$$
A_{ij}= \int_{C_{j}(Q_i)} |T_ib_i||u|, \quad C_{j}(Q_{i}) = 2^{j+1}Q_{i} \setminus 2^jQ_{i}.
$$   
Using the well-known Gaussian upper bounds for the kernels of  
 $tHe^{-tH}$,  ${t>0}$, and 
 $r_i \sim R_i$, we obtain
$$
\|H e^{-t^2 H}  b_i \|_{L^2(C_j(Q_i))} \le \frac{C}{t^{\gamma+2}} e^{-\frac{c 4^jr_i^2}{t^2}}\ \|b_i\|_1
$$
where $\gamma=\frac n 2  $.  
By \eqref{eqcsds3}, 
$\|b_i\|_1 \le
 c
\alpha
R_i |Q_{i}|$,
 hence,
by Minkowski integral inequality, for some appropriate  positive constants $C,c$,
\begin{align*}\| T_ib_i \|_{L^2(C_j(Q_i))} &\le   \int_0^{r_i} \|He^{-t^2 H}  b_i \|_{L^2(C_j(Q_i))} \, dt \\
&\le C\alpha  e^{-c4^j}  |Q_{i}|^{1/  2}.
\end{align*}
Now remark that for any $y \in Q_i$ and any $j\ge 2$, 
$$
\left( \int_{C_j(Q_i)} |u|^{2}\right)^{1/2} \le \left( \int_{2^{j+1}Q_i} |u|^{2}\right)^{1/2} \le  (2^{n(j+1)}|Q_i|)^{1/2}
\big(M(|u|^{2})(y)\big)^{1/2}.
$$
Applying H\"older inequality,     one obtains
$$
A_{ij} \le  C\alpha 2^{nj/2} e^{-c4^j}  |Q_{i}| \big(M(|u|^{2})(y)\big)^{1/2}.
$$
Averaging over $Q_i$ yields
$$
A_{ij} \le C\alpha 2^{nj/2} e^{-c4^j} \int_{Q_i} \big(M(|u|^{2})(y)\big)^{1/2}\, dy.$$
Summing over $j\ge 2$ and $i$, we have
$$
\int| u| \sum_i h_i \le C \alpha \int \sum_i {\bf 1}_{Q_i}(y) \big(M(|u|^{2})(y)\big)^{1/2}\, dy.$$
Using finite overlap \eqref{eqcsds5} of the cubes $Q_i$ and Kolmogorov's inequality, one obtains
$$
\int| u| \sum_i h_i  \le C'N\alpha \big| \cup_i Q_i \big|^{1/2} \||u|^2\|_1^{1/2}.
$$
Hence 
$$
\bigg|\bigg\{x \in \RR^n \setminus \cup_i 4Q_i ; \bigg|\sum_i T_ib_i(x)\bigg| >\frac \alpha 3\bigg\}\bigg| \le C \big| \cup_i Q_i \big| \le 
\frac{C}{\alpha^p}\int |\nabla f|+ |V^{1/2}f|$$
by \eqref{eqcsds5} and \eqref{eqcsds4}.

It remains to handling  the term $B$. Using functional calculus for $H$ one can
compute 
$U_i$ as $r_i^{-1}\psi(r_i^2 H)$ with $\psi$ the holomorphic function on the sector $|\arg z \,| < {\frac \pi 2}$ given
by 
\begin{equation*}\label{eqpsi}
\psi(z)= \int_1^\infty e^{-t^2z} z\, dt.
\end{equation*}
It is easy to show that $|\psi(z)| \le C|z|^{1/2} e^{-c|z|}$, uniformly on subsectors $|\arg z\,  | \le \mu < {\frac \pi 2}$.

Let $q=2$ if $n=1$ and $q=1^*=\frac n {n-1}$ for $n\ge 2$. By Poincar\'e-Sobolev inequality, 
$b_{i}\in L^q$ and 
$$\|b_i\|_q \le c R_i^{1-( n  -\frac n q)}
\|\nabla b_i
\|_1 \le C\alpha R_{i}^{1+ \frac nq}.$$
We invoke the estimate
\begin{equation}\label{eq23}
\left\| \sum_{k\in \ZZ} \psi(4^k H) \beta_k \right\|_q \lesssim \left\|\left(\sum_{k\in \ZZ}  |\beta_k|^2\right)^{1/2} \right\|_q.
\end{equation}
Indeed, by duality, this is equivalent to the Littlewood-Paley inequality
$$
\left\| \left(\sum_{k\in \ZZ} |\psi(4^kH) \beta|^2\right)^{1/2} \right\|_{q'} \lesssim \|\beta\|_{q'}.
$$
For $q=2$, this is a simple estimate using Borel functional calculus on $L^2$ since $H$ is self-adjoint. For $q\ne 2$,  this is  a consequence of the Gaussian estimates for the kernels of $e^{-tH}$, $t>0$ (this was first proved in \cite{ADMc} using the vector-valued version of the work in \cite{DMc}.  See \cite{Aus} for  a more general argument in this spirit or  \cite{LeM} for an abstract proof relying on functional calculus).

To apply \eqref{eq23}, observe that the definitions of $r_i$ and $U_i$ yield 
$$
\sum_i U_ib_i = \sum_{k\in \ZZ} \psi(4^k H) \beta_k
$$
with 
$$\beta_k = \sum_{i, r_i=2^k} \frac{b_i}{r_i}.
$$
Using the bounded overlap property \eqref{eqcsds5}, one has that 
$$
\left\|\left(\sum_{k\in \ZZ}  |\beta_k|^2\right)^{1/2} \right\|_q^q \le C\int \sum_{i}  \frac{|b_i|^q}{r_i^q} .
$$
Using $R_i\sim r_i$, 
$$
\int \sum_{i}  \frac{|b_i|^q}{r_i^q}  \le  C \alpha^q  \sum_{i}|Q_i| .$$
Hence, by  \eqref{eqcsds4}
$$
\bigg|\bigg\{x \in \RR^n; \bigg|\sum_i U_ib_i(x)\bigg| 
>\frac \alpha 3\bigg\}\bigg| \le  C   \sum_{i}|Q_i| \le \frac{C}{\alpha}\int
|\nabla f|+|V^{1/2}f|.$$

We turn to the proof of  \eqref{eq:marcin}.
Fix $1<p<2$ and $f\in C_{0}^\infty(\RR^n)$. 
Choose
$0<\delta<1$ so that $1<p\delta$. Let $\alpha>0$ and apply the Calder\'on-Zygmund decomposition of Lemma \ref{lemmaCZD} to $ f$  with exponent $p\delta$ and threshold $\alpha$.  We may do this since $\|\nabla f\|_{p\delta } + \|V^{1/2}f\|_{p\delta }<\infty$. Of course we do not want to use its value  in a quantitative way.  We obtain that $f=g_{\alpha}+ b_{\alpha}$ with $b_{\alpha}=\sum_{i}b_{i}$.

Write 
\begin{align*}
\|H^{1/2} f\|_{p}^p &= p 2^p \int_{0}^\infty \alpha^{p-1} |\{ x \in \RR^n; |H^{1/2}f(x)| >2\alpha\}| \, d\alpha
\\
&
\le p 2^p \int_{0}^\infty \alpha^{p-1} |\{ x \in \RR^n; |H^{1/2}g_{\alpha}(x)| >\alpha\}| \, d\alpha
\\
&
\qquad + p 2^p \int_{0}^\infty \alpha^{p-1} |\{ x \in \RR^n; |H^{1/2}b_{\alpha}(x)| >\alpha\}| \, d\alpha
\\
&
\le 
I + II
\end{align*}
with
$$
I = C p 2^p \int_{0}^\infty \alpha^{p-1}  \frac {\|\nabla g_{\alpha}\| _2^2 +\|V^{1/2}g_{\alpha}\|_{2}^2} {\alpha^2}   \, d\alpha = I_{g}+I_{v}
$$
and
$$
II= C p 2^p \int_{0}^\infty \alpha^{p-1}  \frac {\|\nabla b_{\alpha}\|_1+ \|V^{1/2}b_{\alpha}\|_{1} } {\alpha}   \, d\alpha= II_{g}+ II_{v},
$$
where $I_{g}$ and $II_{g}$ denote the gradient term in $I$ and $II$ respectively. 
To estimate these integrals, we need to come back to the construction of $ g_{\alpha}$ and $b_{\alpha}$. 
 
Set  $Tf=(|\nabla f|^{p\delta}+  | V^{1/2}f |^{p\delta })^{1/p\delta }$. Write $F_{\alpha}$ as the complement of $\Omega_{\alpha}= \{M(Tf^{p\delta}  ) > \alpha^{p\delta }\}$. Then 
recall that $\nabla g_{\alpha}=  \mathbf{1}_{F_{\alpha}} (\nabla f) +  \mathbf{1}_{\Omega_{\alpha}} h$ where $|h|\le C\alpha$ and $|\nabla f |\le \alpha$ on $F_{\alpha}$. Thus $I_{g}$ splits into $I_{g1}+I_{g2}$ according to this decomposition. The treatment
of $I_{g1}$ is done using the definition of $F_{\alpha}$, Fubini's theorem and $p<2$ as follows:  
\begin{align*}
I_{g1}&= \frac{C p 2^p}{2-p} \int    |\nabla f|^2 \ \left({ M}( Tf^{p\delta})\right)^{\frac {p-2}{p\delta }}   \
\\
&
\le \frac{C p 2^p}{2-p} \int    |\nabla f|^p ,
\end{align*}
where we used 
$$|\nabla f|^2= |\nabla f|^p\ |\nabla f|^{2-p} \le |\nabla f|^p\ (Tf^{p\delta })^{\frac{2-p}{p\delta} }  \le  |\nabla f|^p\ \left({M}( Tf^{p\delta})\right)^{\frac {2-p}{p\delta }} $$ almost everywhere. 
For $I_{g2}$, we  use the bound of $h$ to obtain
\begin{align*}
I_{g2}&\le  C p 2^p \int_{0}^\infty \alpha^{p-1}  |\Omega_{\alpha}|   \, d\alpha
\\
&
= C 2^p \int \left(M( Tf^{p\delta})\right)^{\frac {1}{\delta }}   
\\
&
\le C \int    |\nabla f|^{p}+  | V^{1/2}f |^{p } 
\end{align*}
by the strong type $(\frac 1 \delta , \frac 1 \delta )$ of the maximal operator.

Next, we turn to the term $II_{g}$. We have $\nabla b_{\alpha}= \mathbf{1}_{\Omega_{\alpha}}(\nabla f)  -  \mathbf{1}_{\Omega_{\alpha}} h  $ so that $II_{g}  \le (II_{g1}+I_{g2})$ and $I_{g2}$ is already controlled. For $II_{g1}$ we have by using H\"older's inequality 
and  the strong type $(\frac 1 \delta , \frac 1 \delta )$ of the maximal operator
\begin{align*}
II_{g1}&= \frac{C p 2^p}{p-1} \int     |\nabla f| \ \left({M}( Tf^{p\delta})\right)^{\frac {p-1}{p\delta }}   
\\
&
\le \frac{C p 2^p}{p-1} \left(  \int    |\nabla f|^p  \right)^{1/p} \left( \int   \left( {M}( Tf ^{p\delta})\right)^{(\frac {p-1}{p\delta }) p'}   \right)^{1/p'}
\\
&
\le C \int    |\nabla f|^{p}+  | V^{1/2}f |^{p}.
\end{align*}

It remains to look at $I_{v}$ and $II_{v}$. Recall that $g_{\alpha}= f$ on $F_{\alpha}$ and $g_{\alpha}=h_{\alpha}$ on $\Omega_{\alpha}$, and we have proved
$\int V|h_{\alpha}|^2 \le C\alpha^2 |\Omega_{\alpha}|$. Hence, $I_{v}$ splits as $I_{v1}+I_{v2}$. First, 
\begin{align*}
I_{v1}&= \frac{C p 2^p}{2-p} \int    |V^{1/2} f|^2 \ \left({ M}( Tf^{p\delta})\right)^{\frac {p-2}{p\delta }}   \
\\
&
\le \frac{C p 2^p}{2-p} \int    |V^{1/2} f|^p.
\end{align*} 
with the similar argument as for $I_{g1}$. Next, 
\begin{align*}
I_{v2}&\le  C p 2^p \int_{0}^\infty \alpha^{p-1}  |\Omega_{\alpha}|   \, d\alpha
\\
&
= C 2^p \int \left(M( Tf^{p\delta})\right)^{\frac {1}{\delta }}   
\\
&
\le C \int    |\nabla f|^{p}+  | V^{1/2}f |^{p }. 
\end{align*}

Now, $b_{\alpha}= f-g_{\alpha}= f - h_{\alpha}$ on $\Omega_{\alpha}$ and $b_{\alpha}=0$ on $F_{\alpha}$. Hence, $II_{v}\le II_{v1}+I_{v2}$ and
\begin{align*}
II_{v1}&= \frac{C p 2^p}{p-1} \int     |V^{1/2} f| \ \left({M}( Tf^{p\delta})\right)^{\frac {p-1}{p\delta }}   
\\
&
\le \frac{C p 2^p}{p-1} \left(  \int    |V^{1/2} f|^p  \right)^{1/p} \left( \int   \left( {M}( Tf ^{p\delta})\right)^{(\frac {p-1}{p\delta }) p'}   \right)^{1/p'}
\\
&
\le C \int    |\nabla f|^{p}+  | V^{1/2}f |^{p}.
\end{align*}
 This concludes the proof of item 3 of Theorem \ref{th:2}.

\section{Estimates for weak solutions}\label{sec:weaksol}

  In this section,  $Q$ denotes a cube, $R$ its radius, and $u$  a  weak solution of $-\Delta u + Vu=0$ in a neighborhood of $\overline{2Q}$. Recall that under the assumption
  $V\ge 0$, we have the mean value inequality
  \begin{equation}
\label{eq:mvi}
\sup_{Q} |u| \le C(r,n,\mu) \big(\av_{\mu Q} |u|^r\big)^{1/r}
\end{equation}
for any $0<r<\infty$ and $1<\mu\le 2$.  And we have also shown a mean value inequality
against arbitrary $A_{\infty}$ weights.

 We state some further estimates  that are interesting in their own right assuming $
 V\in A_{\infty}$.  By splitting real and imaginary parts,  we may  suppose $u$ real-valued. All constants are  independent of $Q$ and $u$ but they may depend on $V$ through the constants in the $A_{\infty}$ condition or the $B_{q}$ condition when assumed. 

\begin{lem}
\label{lem:1}
For all $1 \le \mu < \mu' \le 2$ and $k>0$, there is a constant $C$ such that 
$$ 
\av_{\mu Q}|u|^2 \le \frac {C} { (1+ R^2\av_{Q}V)^k} \big(\av_{\mu' Q} |u|^2\big).
$$
and
$$ 
\av_{\mu Q} (|\nabla u|^2 + V|u|^2) \le \frac {C} { (1+ R^2\av_{Q}V)^k} \big(\av_{\mu' Q} (|\nabla u|^2 + V|u|^2) \big).
$$
\end{lem}

\begin{lem}
\label{lem:2}
For all $1<\mu \le 2$ and $k>0$,  there is a constant $C$ such that 
$$
(R\av_{Q}V)^2 \, \av_{Q} |u|^2 \le \frac {C} { (1+ R^2\av_{Q}V)^k}  \big(\av_{\mu Q}  ( V|u|^2 )\big).
$$ 
\end{lem}

\begin{lem}
\label{lem:3}
For all $1<\mu\le 2$, $k>0$  and $\sup(n,2) <p<\infty$, there is a constant $C$ such that 
$$
(R\av_{Q}V)^2 \, \av_{Q} |u|^2 \le  \frac {C} { (1+ R^2\av_{Q}V)^k}  \big( \av_{\mu Q}   |\nabla u|^p \big)^{2/p}.
$$
\end{lem}

\begin{lem}
\label{lem:4}  
Assume $V \in B_{q}$  and set $\tilde q = \inf (q^*, 2q)$. For all $1< \mu\le 2$ and $k>0$  there is a constant $C$ such that
$$
\big( \av_{ Q}   |\nabla u|^{\tilde q} \big)^{1/\tilde q} \le   \frac {C} { (1+ R^2\av_{Q}V)^k} \, \big(\av_{\mu Q}  ( |\nabla u|^2 + V|u|^2 )\big)^{1/2}. 
$$
\end{lem}

\begin{lem}
\label{lem:5}  
Assume $V \in B_{q}$. For all $1< \mu\le 2$, if  $n/2 \le q <n$ then  there is a constant $C$ such that
$$
\big( \av_{ Q}   |\nabla u|^{ q^*} \big)^{1/ q^*} \le  C \, \big(\av_{\mu Q}   |\nabla u|^2\big)^{1/2}, 
$$
and if  $q \ge n$ then there is a constant $C$ such that
  $$
 \sup_{ Q}   |\nabla u|  \le  C \, \big(\av_{\mu Q}   |\nabla u|^2\big)^{1/2}. 
$$
\end{lem}

\begin{lem}
\label{lem:6}
Assume $V \in B_{q}$. For all  $1< \mu\le 2$ and $k>0$, if  $n/2 \le q <n$ then there is a constant $C$ such that
$$
\big( \av_{ Q}   |\nabla u|^{ q^*} \big)^{1/ q^*} \le  \frac {C} { R(1+ R^2\av_{Q}V)^k} \big( \sup_{\mu Q} |u|\big). 
$$
and if  $q \ge n$ then  there is a constant $C$ such that
$$
 \sup_{ Q}   |\nabla u|  \le    \frac {C} { R(1+ R^2\av_{Q}V)^k}  \big(\sup_{\mu Q} |u|\big). 
$$

\end{lem}

\begin{lem}
\label{lem:7} 
Assume $V \in B_{q}$ with $q>1$ and $ q\ge n/2$. For all  $1< \mu\le 2$ and $k>0$ there is a constant $C$ such that
$$
(R  \av_{ Q} V)^2  \,  \av_{ Q} | u|^2 \le   \frac {C} { (1+ R^2\av_{Q}V)^k}  \big(\av_{\mu Q} |\nabla u|^2\big). 
$$
\end{lem}

\begin{lem}
\label{lem:8} 
In Lemma \ref{lem:5}, the constant $C$ can be replaced by $C(1+R^2\av_{Q}V)^{-k}$ for any $k>0$.
\end{lem}

Let us postpone the proofs and make some remarks concerning these inequalities.
\begin{rem}\label{rem:est}
1) Lemma \ref{lem:5} is a weak reverse H\"older inequality for  the gradient of weak solutions.  It improves over Lemma  \ref{lem:4} in the fact that the right hand side does not have terms involving $V|u|^2$ but this is under the assumption $q \ge n/2$.  
Using self-improvement of weak reverse H\"older inequalities $($see \cite[Theorem 2]{IN}$)$, we may replace the exponent $2$ in the right hand sides by any $0<p<2$.

2) We do not know if Lemma \ref{lem:5} holds for $q<n/2$. 

3) In Lemma \ref{lem:4},   note that $\tilde q= q^*<2q$ when $q<n/2$ and it  would be natural  the estimate holds  for the larger exponent $2q$. 

4) Lemma \ref{lem:7} is a Poincar\'e type inequality for weak solutions. As $\sup_{Q}  | u|$ can be compared to $\big(\av_{Q} | u|^2\big)^{1/2}$, we see that it is a converse to the Caccioppoli inequality in the regime $R^2\av_{Q}V\ge 1$. 

5)   Except for Lemma \ref{lem:1} and  \ref{lem:6} which  are closely related to  Lemma 4.6 and Remark 4.9 in \cite{Sh1}, these lemmata  appear to be new.

\end{rem}

\paragraph{\bf Proof of Lemma \ref{lem:1}}  There is nothing to prove if $R^2\av_{Q}V \le 1$ and we assume $R^2 \av_{Q}V \ge 1$. 
The well-known Caccioppoli type argument yields for $1\le \mu <\mu' \le 2$
\begin{equation}
\label{eq:caccio}
\int_{\mu Q} |\nabla u|^2 + V|u|^2 \le \frac {C} {R^2} \int_{\mu' Q} |u|^2.
\end{equation}
The improved Fefferman-Phong inequality of Lemma \ref{lem:FP} and the fact that the averages of $V$ on $\mu Q$ with $1\le \mu \le 2$ are all uniformly comparable tell us  for  some $\beta>0$,   
$$
 \frac{ 1}{ R^{2}}  \int_{\mu Q} |u|^2 \le \frac C {(R^2 \av_{Q}V)^\beta}\int_{\mu Q}  |\nabla u|^2 + V|u|^2  .
$$ 
The desired estimates follow readily by iterating  these two inequalities.

\

\paragraph{\bf Proof of Lemma \ref{lem:2}} Using Lemma \ref{lem:1} with $k>1$ and $1<\mu' <\mu$ and then Lemma \ref{lem:subw}, we have, 
$$
(R \av_{Q} V)^2 \av_{ Q}|u|^2 \le \frac {C  \av_{Q}V\,\av_{ \mu' Q}|u|^2} { (1+ R^2\av_{Q}V)^{k-1}}  \le \frac {C  \av_{\mu' Q}V\,\sup_{ \mu' Q}|u|^2} { (1+ R^2\av_{Q}V)^{k-1}}   \le \frac {C  \av_{\mu Q} (V|u|^2)} { (1+ R^2\av_{Q}V)^{k-1}}  .
$$
 
 \

\paragraph{\bf Proof of Lemma \ref{lem:3}} Of course, if $\av_{\mu Q} |\nabla u|^p
=\infty$ there is nothing to prove. Assume, therefore, that $\av_{\mu Q} |\nabla u|^p
<\infty$.  
Let $1<\nu <\mu$ and $\eta$ be a smooth non-negative function, bounded by 1, equal to 1 on $\nu Q$ with support on $\mu Q$ and whose gradient is bounded by $ C/R$ and Laplacian by $C/R^2$.
Integrating the equation $-\Delta u + Vu =0$ against $u\eta^2$, we find
$$
\int  |\nabla u|^2 \eta^2 + V|u|^2 \eta^2 = 2 \int \nabla u \cdot \nabla \eta\, u \eta 
\le \frac C R \bigg( \int_{\mu Q} |\nabla u|^2 \bigg)^{1/2} \bigg( \int |u|^2 \eta^2\bigg)^{1/2},
$$
hence
$$
X \le  {C\, {(R^2\av_{Q}V)^{1/2}}  |\mu Q|^{1/2}\, Y^{1/2}\, Z^{1/2}} $$
where we have set $X=(R^2\av_{Q}V) \int V |u|^2 \eta^2$,  $Y= \big(\av_{\mu Q} |\nabla u|^p\big)^{2/p}$
and $Z= \av_{Q}V \int |u|^2 \eta^2$.  By Morrey's embedding theorem, $u$ is H\"older continuous with exponent $\alpha= 1- n/p$ and for all $x,y\in \mu Q$, 
$$
|u(x) -u(y)| \le C \bigg( \frac{|x-y|}R\bigg)^\alpha \, R\, \big(\av_{\mu Q} |\nabla u|^p\big)^{1/p}= C \bigg( \frac{|x-y|}R\bigg)^\alpha \, R\, Y^{1/2}.
$$ We pick $y\in \overline Q$ such that $|u(y)|= \inf_{Q}|u|$. Then
\begin{align*}
Z=\av_{Q}V \int |u|^2 \eta^2 &\le 2  (\av_{Q} V) \inf_{Q}|u|^2 \int \eta^2 + 2 (\av_{Q} V) \int |u(x) -u(y)|^2 \eta^2(x)\, dx
\\
& \le 2 \big(\av_{Q} (V|u|^2)\big) \int \eta^2 + C (\av_{Q} V) R^2 Y \, \int 
 \bigg( \frac{|x-y|}R\bigg)^{2\alpha} \eta^2(x)\, dx
 \\
 & \le C   \big(\av_{Q} (V|u|^2)\big)  |Q| + C (\av_{Q} V) R^2 Y\,|\mu Q|  \\
 & \le C  \int V|u|^2\eta^2 + C (\av_{Q} V) R^2 Y\,|\mu Q|,
 \end{align*}
 where, in the penultimate inequality, we used the support condition on $\eta$ and $0\le \eta \le 1$, and in the last,    $\eta=1$ on $ Q$. 
 Using the previous inequalities, we obtain
$$X  \le C  |\mu Q| ^{1/2} \, Y^{1/2} \, \big( CX + C(R^2\av_{Q}V)^2  |\mu Q| Y\big)^{1/2}
 $$
 which, by $2ab \le \ep^{-1} a^2 + \ep b^2$ for all $a,b \ge 0$ and $\ep>0$, implies
 $$
 X \le C(1 + R^2 \av_{Q}V )^2 \,  |\mu Q| \, Y.
 $$
 Next, let $1< \nu'<\nu$. Using $\eta=1$ on $\nu Q$,  Lemma \ref{lem:subw} and Lemma \ref{lem:1},
 $$
 \int  V|u|^2 \eta^2 \ge \int_{\nu Q} V |u|^2 \ge C \av_{\nu' Q } V \, \int_{\nu' Q} |u|^2 \ge C   (\av_{Q } V)  (1 + R^2 \av_{Q}V )^k \int_{ Q} |u|^2,
 $$
 hence
 $$
 X \ge C (R \av_{Q}V)^2  (1 + R^2 \av_{Q}V )^k \int_{Q} |u|^2.
 $$
 The upper and lower bounds for $X$ yield the lemma. 
 
 \
 
\paragraph{\bf Proof of Lemma \ref{lem:4}}  First note that if $q\le \frac{2n}{n+2}$ then $\tilde q\le 2$ and the conclusion (useless for us) follows by a mere H\"older inequality. Henceforth, we assume $q>\frac{2n}{n+2}$. Also, by Lemma \ref{lem:1}, it suffices to obtain the estimate with $k=0$. Let us assume $\mu=2$ for simplicity of the argument.
Let $v$ be the harmonic function  on $ 2 Q$ with $v=u$ on $\partial( 2 Q)$ and set $w=u-v$ on $ 2 Q$. Since $w=0$ on $\partial( 2 Q)$, we have 
$$
 (\av_{ 2 Q}   |\nabla w|^2 \big)^{1/2}  \le 
(\av_{ 2 Q}   |\nabla u|^2 \big)^{1/2} .
$$
By elliptic estimates for harmonic functions, we have for all $2\le p \le \infty$, and in particular for $p=\tilde q$,  
$$
\big ( \av_{Q} |\nabla v|^p \big)^{1/p} \le C (\av_{ 2 Q}   |\nabla v|^2 \big)^{1/2}  \le 2 C 
(\av_{ 2 Q}   |\nabla u|^2 \big)^{1/2}.
$$

 Let $1<\mu <2$ and  $\eta$ be a smooth non-negative function, bounded by 1, equal to 1 on $ Q$ with support contained in $\mu Q$ and whose gradient is bounded by $ C/R$ and Laplacian by $C/R^2$.  As $\Delta w= \Delta u = V u $ on $2Q$, we have
$$
\Delta(w\eta) =  Vu \eta +  2 \nabla w\cdot \nabla \eta  + w \Delta\eta  \quad {\rm on \ } \RR^n.
$$
Hence, if $n\ge 2$ by Green's representation for the Laplace equation 
$$
\nabla (w\eta)(x) = \int_{\RR^n} \nabla\Gamma(x-y) \big[(Vu\eta)(y) + 2\nabla w (y) \cdot \nabla \eta(y)  + w(y) \Delta\eta(y)\big]\, dy= I+II+III
$$
where $\Gamma$ is the fundamental solution of $\Delta$ so  that $|\nabla\Gamma(x) | \le C |x|^{1-n}$.   Since $\tilde q \le q^*$, we have
$$
\big( \av_{Q} |\nabla w|^{\tilde q} \big)^{1/\tilde q } \le \big( \av_{Q} |\nabla w|^{ q^*} \big)^{1/ q^* }
$$
so that it suffices to bound the latter integral.
Using support conditions on $\eta$, we obtain the pointwise bounds for $x\in Q$,
$$
II \le  C \av_{ 2 Q} |\nabla w| \le C \big( \av_{ 2 Q} |\nabla w|^2 \big)^{1/2} \le C \big( \av_{ 2 Q} |\nabla u|^2 \big)^{1/2}
$$
and 
$$
III \le \frac C R \av_{ 2 Q} | w| \le C \big( \av_{ 2 Q} |\nabla w|^2 \big)^{1/2} \le C \big( \av_{ 2 Q} |\nabla u|^2 \big)^{1/2}
 $$
 where we used Poincar\'e inequality for $w$ on $2Q$ as $w=0$ on the  boundary. 
By the $L^q-L^{q^*}$ boundedness of the Riesz potential
$$
\bigg(\int_{\RR^n} I^{q^*}\bigg)^{1/q^*} \le C\bigg( \int_{
\RR^n} |Vu\eta|^q \bigg)^{1/q}  \le C\bigg( \int_{
\mu Q} |V|^q \bigg)^{1/q} \sup_{\mu Q } |u|.
$$
Normalizing by taking averages and using the $B_{q}$ condition on $V$ yields
\begin{equation}
\label{eq:I}
\big (\av_{Q} I^{q^*}\big)^{1/q^*} \le C R \, \av_{\mu Q} V \sup_{\mu Q } |u|.
\end{equation}
Now, if $\mu < \mu' < 2$,  subharmonicity of $|u|^2$ and Lemma \ref{lem:subw} yield
$$
R\, \av_{\mu Q} V \sup_{\mu Q}|u| \le  C R\,  \av_{\mu' Q} V  \,  \big(\av_{\mu' Q} |u|^2 \big)^{1/2}
$$
 which by Lemma \ref{lem:2} is bounded by 
 $
C  \big(  \av_{2Q} ( V|u|^2)\big)^{1/2}.
$
Gathering the estimates obtained for $\nabla v$ and $\nabla w$,  the lemma is proved when $n\ge 2$.

When $n=1$,  we have 
$$
(w\eta)'(x) = - \int_{x}^\infty Vu\eta + 2 w'\eta' + w \eta''
$$
and we obtain readily for $x\in Q$, 
$$
|w'(x)| \le C R (\av_{\mu Q} V) \sup_{\mu Q} |u| + C \big( \av_{\mu Q}|w'|^2\big)^{1/2}.
$$
The rest of the proof is as before. 

\

\paragraph{\bf Proof of Lemma \ref{lem:5}} Assume $n/2 < q < n$. The previous lemma shows that $ \av_{\mu' Q} |\nabla u|^{ \tilde q}  <\infty$ for all $1<\mu'\le \mu$. As $\tilde q =2q>n$, Lemma \ref{lem:3} applies and using it  with $k=0$ instead of Lemma \ref{lem:2} in the previous argument, we obtain, 
$$
\big( \av_{Q} |\nabla w|^{ q^*} \big)^{1/ q^* } \le C \big( \av_{\mu Q} |\nabla u|^{ 2q} \big)^{1/ 2q }.
$$
As the similar estimate holds for $v$ in place of $w$, we obtain
$$
\big( \av_{Q} |\nabla u|^{ q^*} \big)^{1/ q^* } \le C \big( \av_{\mu Q} |\nabla u|^{ 2q} \big)^{1/ 2q }.
$$
Note that this inequality holds not just for $Q$ but for all cubes $Q'$ with $\overline{2Q'}$ contained in the open set where $u$ is a weak solution. As $q^*>2q$, this set of inequalities self-improves   with $2q$ replaced by any $0<p<2q$ (see \cite{IN})  and, in particular, 
$$
\big( \av_{Q} |\nabla u|^{ q^*} \big)^{1/ q^* } \le C \big( \av_{\mu Q} |\nabla u|^{ 2} \big)^{1/ 2 }.
$$
If $q\ge n$ and $n\ge 2$, then  we may as well consider $q>n$.  Then \eqref{eq:I} becomes
$$\sup_{Q} I \le C R \, \av_{\mu Q} V \sup_{\mu Q } |u|
$$
so that the pointwise bound for $\nabla u$ follows by Lemma \ref{lem:3}.
If $n=1$, we already obtained a pointwise bound for $\nabla u$ and again Lemma \ref{lem:3} applies. 

\

\paragraph{\bf Proof of Lemma \ref{lem:6}} It suffices to incorporate the  Caccioppoli inequality \eqref{eq:caccio} in the inequalities of Lemma \ref{lem:6}.

\

\paragraph{\bf Proof of Lemma \ref{lem:7}}  It suffices to combine Lemma \ref{lem:3} and Lemma \ref{lem:5}.

\

\paragraph{\bf Proof of Lemma \ref{lem:8}} It suffices to see the case $R^2\av_{Q }V \ge 1$. Then, combine Lemma \ref{lem:6}, the mean value inequality \eqref{eq:mvi}  with $r=2$ and Lemma \ref{lem:7}. 

\section{Riesz transforms}

This section is concerned with the proof of  Theorem \ref{th:2}, item 3. We present an argument inspired by \cite{Sh2} which also gives us  a second proof of part of  item 1\footnote{In this section, $L^p$ denotes either $L^p(\RR^n, \CC)$ or $L^p(\RR^n, \CC^n)$.}.

\subsection{A reduction}

We know that it suffices to establish the boundedness of $\nabla H^{-1/2}$ and of $V^{1/2}H^{-1/2}$ on $L^p$ for the appropriate ranges of $p$. As already observed, the case $1<p\le 2$ is already taken care of with no assumption on $V$. We henceforth assume $p>2$ and $V\in A_{\infty}$.

 By duality, we know that $H^{-1/2}\div$ and $H^{-1/2}V^{1/2}$ are bounded on $L^p$. Thus, if $\nabla H^{-1/2}$ is also bounded on $L^p$, we obtain that $\nabla H^{-1} \div$ and $\nabla H^{-1}V^{1/2}$ are bounded on $L^p$. 

Reciprocally, if $\nabla H^{-1} \div$ and $\nabla H^{-1}V^{1/2}$ are bounded on $L^p$, then their adjoints are bounded on $L^{p'}$. Thus, if $F\in C_{0}^\infty(\RR^n, \CC^n)$, 
\begin{align*}
\|H^{-1/2} \div F\|_{p'} &= \|H^{1/2} H^{-1}\div F\|_{p'}\\
&\le C (\|\nabla H^{-1} \div F\|_{p'}
+ \| V^{1/2} H^{-1} \div F\|_{p'}) \le C \|F\|_{p'}
\end{align*}
where the first inequality follows from  item 2 of Theorem \ref{th:2}. Hence, by duality, $\nabla H^{-1/2}$ is bounded on $L^p$. 

The same treatment can be done on $V^{1/2} H^{-1/2}$.  We have obtained 
\begin{lem} If $V\in A_{\infty}$ and $p>2$,  the $L^p$ boundedness of  $\nabla H^{-1/2}$  is equivalent to that of $\nabla H^{-1} \div$ and $\nabla H^{-1}V^{1/2}$, and  the $L^p$ boundedness of  $V^{1/2} H^{-1/2}$  is equivalent to that  of  $V^{1/2} H^{-1} V^{1/2}$ and $V^{1/2} H^{-1} \div$.
\end{lem}

 It suffices therefore to establish part of Corollary \ref{cor:3} namely,

 \begin{pro} \label{thC'} Assume that $V \in B_q$ for some $q>1$. Then for 
$2<p\le 2(q+\ep)$, for some $\ep>0$ depending only on $V$,  $f\in C^\infty_{0}(\RR^n,\CC)$ and $F\in C_{0}^\infty(\RR^n, \CC^n)$, 
$$
\|V^{1/2} H^{-1} V^{1/2}f\|_p  \le C_{p} \|f\|_{p}, \quad  \|V^{1/2} H^{-1} \div F\|_p \le C_p\|F\|_p.
$$
\end{pro}

\begin{pro} \label{thD'}
Assume that $V \in B_q$  for some $q>1$. 
Then for $2<p\le  q^*+\ep$ for some $\ep>0$ depending only on $V$, $f\in C^\infty_{0}(\RR^n,\CC)$ and $F\in C_{0}^\infty(\RR^n, \CC^n)$, 
$$
\|\nabla  H^{-1} V^{1/2}f\|_p  \le C_{p} \|f\|_{p}, \quad  \|\nabla H^{-1} \div F\|_p \le C_p\|F\|_p.
$$
\end{pro}

  The interest of such a reduction is that this allows us to use properties of weak solutions of  $H$.
  
   Note that Proposition \ref{thD'} is void if $q\le \frac {2n}{n+2} $ as $q^*\le 2$. Note also that $q^*<2q$ exactly when $q<n/2$. In  this case, this statement yields a smaller range  than 
  the interpolation method in Section \ref{sec:interpolation}. 
 
\subsection{Proof of Proposition \ref{thC'}}

%

%

 Fix a cube $Q$ and and let $f\in C^\infty_{0}(\RR^n)$ supported away from $4Q$. Then $u=H^{-1}V^{1/2}f$ is well defined on $\RR^n$ with 
$\|V^{1/2} u\|_{2}+ \|\nabla u \|_{2} \le \|f\|_{2}$ by construction of $H$ and 
$$
\int_{\RR^n} V u \varphi + \nabla u \cdot \nabla \varphi = \int_{\RR^n} V^{1/2} f \varphi
$$ for all $\varphi\in L^2$ with $\|V^{1/2} \varphi\|_{2}+ \|\nabla \varphi \|_{2}<\infty$.
In particular, the support condition on $f$ implies that $u$ is a weak solution of $-\Delta u + Vu=0$ in $4Q$, hence $|u|^2$ is subharmonic on $4Q$.   Let $r$ such that $V\in B_{r}$ and note that $V^{1/2} \in B_{2r}$ (see section \ref{sec:Ainfty}). By Lemma  \ref{lem:rh2q} with $w=V^{1/2}$ $f=|u|^2$ and $s=1/2$, we have
$$
\big(\av_{Q} (V^{1/2}|u| )^{2r} \big)^{1/2r} \le C \, \av_{\mu Q} (V^{1/2} |u|).
$$
Thus,  \eqref{T:shen} holds with $T=V^{1/2}H^{-1} V^{1/2}$, $q_{0}=2r$, $p_{0}=2$ and $S=0$. By Theorem \ref{theor:shen}, 
$V^{1/2}H^{-1} V^{1/2}$ is bounded on $L^p$ for $2<p<2r$. 

The argument is the same  for $V^{1/2}H^{-1}\div$. This finishes the proof.

\subsection{Proof of Proposition \ref{thD'}}

 We assume $q>\frac{2n}{n+2}$, that is $q^*>2$. otherwise there is nothing to prove. We consider first the operator $\nabla H^{-1} V^{1/2}$.

  Assume $q<n/2$. 
Fix a cube $Q$ and and let $f\in C^\infty_{0}(\RR^n)$  supported away from $4Q$. Let $u=H^{-1}V^{1/2}f$. As before,  the support condition on $f$ implies that $u$ is a  weak solution of $-\Delta u + Vu=0$ in $4Q$.  Thanks to Lemma \ref{lem:4}, 
\eqref{T:shen} holds with $T=\nabla H^{-1} V^{1/2}$, $q_{0}= q^*$ and $S=V^{1/2}H^{-1}V^{1/2}$.
As $S$ is bounded on $L^{q^*}$ by Proposition \ref{thC'} and $2<q^*\le 2 q$, Theorem \ref{theor:shen} implies  that 
$\nabla H^{-1} V^{1/2}$ is bounded on $L^p$ for $2<p <q^*$. Finally, by the self-improvement of reverse H\"older estimates we can replace $q$ by a slightly larger value and, therefore, $L^p$ boundedness  for $p< q^*+\ep$ holds.

Assume next that $ n/2 \le q < n$. in this case, $q^*\ge 2q$. Again, we may as well assume $q>n/2$. Then,  Lemma \ref{lem:5}  yields, this time, \eqref{T:shen}  with $T=\nabla H^{-1} V^{1/2}$, $q_{0}= q^*$ and $S=0$. Hence, Theorem \ref{theor:shen} implies  that 
$\nabla H^{-1} V^{1/2}$ is bounded on $L^p$ for $2<p< q^*$.  Again, by self-improvment
of the $B_{q}$ condition, it holds for  $p<q^*+\ep$.

Finally, if $q\ge n$, then,  Lemma \ref{lem:5}  yields  \eqref{T:shen} for any $2<q_{0}<\infty$   with $T=\nabla H^{-1} V^{1/2}$ and $S=0$. Hence, Theorem \ref{theor:shen} implies  that 
$\nabla H^{-1} V^{1/2}$ is bounded on $L^p$ for $2<p < \infty$. 

The argument is the same  for $\nabla H^{-1}\div$ and this finishes the proof.


\section{$L^p$ Domains of $H$ and $H^{1/2}$}

\paragraph{\bf Proof of Corollary \ref{cor:1}} It is known that $-\Delta+V$ defined on $\coinf(\RR^n)$ is essentially m-accretive on $L^p(\RR^n)$ if $V\in L^p_{loc}(\RR^n)$.
The domain of its extension is $\{u\in L^p(\RR^n)\, ; \, -\Delta u + Vu \in L^p(\RR^n)\}$ with norm $\|u\|_{p}+ \|-\Delta u + Vu \|_{p}$. By \eqref{eq:Lp} this norm is equivalent to 
$\|u\|_{p}+ \|\Delta u\|_{p} + \|Vu \|_{p}$ on $\coinf(\RR^n)$ when $V\in B_{p}$. The result follows.

\

\paragraph{\bf Proof of Corollary \ref{cor:2} }

Let $E^{p}(\RR^n)= \D_{p}(\nabla)\cap \D_{p}(V^{1/2}) = W^{1,p}(\RR^n) \cap L^p(\RR^n,V^{p/2})$.
Let us begin with the following lemma.

\begin{lem} If $1< p<\infty$ and $V^{p/2} \in L^1_{loc}(\RR^n)$, then 
$C_{0}^\infty(\RR^n)$ is dense in $ E^{p}(\RR^n)$. \end{lem}  

Indeed, for $p=2$ this is a well-known fact as $C_{0}^\infty(\RR^n)$  is a core of the form 
domain of $-\Delta+V$. The proof of this fact (see, for instance, \cite[pp. 157-158]{Da3}) adapts to any  $p$ with $1<p<\infty$. 

We also remark that under the assumption $V\in L^{1}_{loc}$,  $-\Delta+V$ has a bounded holomorphic functional calculus on $L^p(\RR^n)$  for $1<p<\infty$ (\cite{DR}), and in particular,  $\|(-\Delta + V+1)^{1/2} u \|_{p} \sim \|(-\Delta + V)^{1/2} u \|_{p} + \|u\|_{p}$ for all $u \in C_{0}^\infty(\RR^n)$. Thus, it suffices to find the domain of $(-\Delta + V+1)^{1/2}$.

Now, assume $V\in A_{\infty}$ and $1<p<2$ or $V\in B_{p/2}$ and $2<p<\infty$. We have shown that $\|(-\Delta + V)^{1/2} u \|_{p} \sim \|\nabla u \|_{p}+ \|V^{1/2}u\|_{p}$ for $u\in C_{0}^\infty(\RR^n)$.   Thus,  using this and the lemma, $(-\Delta + V+1)^{1/2}$ has a bounded extension from $E^{p}(\RR^n)$ to $L^p(\RR^n)$ and this extension is invertible. This proves the result. 

\begin{rem}
It is not hard to show that the $L^p$-domain  $(1<p<\infty)$ of $(-\Delta+V)^{1/2}$ coincides with the domain of the square root of $($minus$)$ the infinitesimal generator of the semigroup $(e^{-tH})_{t>0}$ seen as an analytic and $C_{0}$-semigroup on $L^p$. 
\end{rem}

\section{Some facts about $A_{\infty}$ weights}\label{sec:Ainfty}

That $V\in A_{\infty}$ implies $V^{s} \in B_{1/s}$ for $0<s<1$ was first observed  implicitely in \cite{ST}. We give a direct proof for convenience.

\begin{pro} Let $V$ be a nonnegative measurable function. Then the followings are equivalent:
\begin{enumerate}
  \item $V\in A_{\infty}$.
  \item For all $s\in (0,1)$, $V^s\in B_{1/s}$.
  \item There exists $s\in (0,1)$, $V^s\in B_{1/s}$.
\end{enumerate}

\end{pro}

\paragraph{\bf Proof} If $V^s\in B_{1/s}$ for some $s\in (0,1)$, then by the self-improvement property of the $B_{q}$
class, $V^s \in B_{\ep+1/s}$ for some $\ep>0$. Hence, $V\in B_{1+s\ep}$, which implies
$V\in A_{\infty}$. Thus, (2) implies (3) implies (1).

Assume $V\in A_{\infty}$ and $s\in (0,1)$. 
Since $A_{\infty}$ weights satisfy a reverse H\"older inequality, there is 
$r>1$ such that $V\in B_{r}$. Hence, for $A>1$ and any cube $Q$, the set $E_{Q}=\{x\in Q\, ; \, 
V^s(x) > A \av_{Q } V^s\} $ satisfies  
$$
\frac{\int_{E_{Q}} V} {\int_{Q}V} \le C \left(\frac {|E_{Q}|}{|Q|}\right)^{1/r'}.
$$
Since ${|E_{Q}|} \le   A^{-1}  |Q|$, we obtain  
$\int_{E_{Q}} V \le C A^{-1/r'} \int_{Q} V $.  

Choose $A$ such that  $CA^{-1/r'} \le 1/2$. We have
$$\int_{Q} V = \int_{Q\setminus E_{Q}} V+ \int_{E_{Q}} V   \le (A \av_{Q} V^s)^{1/s} |Q|+ \frac 1 2  \int_{Q} V$$
which yields
$$
\av_{Q} V \le 2 (A \av_{Q} V^s)^{1/s}.
$$

\end{document}